\documentclass[11,5pt, a4paper]{article}

\usepackage[utf8]{inputenc}
\usepackage{titlesec}
\usepackage[english]{babel}
\usepackage{graphicx}
\usepackage{float}
\usepackage{fancyhdr}
\usepackage[matrix,arrow,curve]{xy}
\usepackage{pstricks} 
\usepackage{amsmath,amsfonts,verbatim,afterpage,theorem,euscript,mathrsfs,amssymb}
\usepackage{amsfonts}
\usepackage{amssymb}
\usepackage{array}
\usepackage{dsfont}
\usepackage{hyperref}
\usepackage{authblk}

\newtheorem{Proposition}{Proposition}[section]

\newtheorem{Lemma}{Lemma}[section]
\newtheorem{Theorem}{Theorem}

\newtheorem{Remark}{Remark}[section]

\def \vu{\textbf{u}}

\def \vv{\textbf{v}}

\newcommand\vN{{\bf  \nabla}}

\def \P{\mathbb{P}}

\def \Rt{\mathbb{R}^{3}}
\def \R{\mathbb{R}}

\def \ve{{\bf e}}
\def \finpv{\hfill $\blacksquare$  \newline }
\def \pv{{\bf{Proof.}}~}
\def \ds{\displaystyle}

\begin{document}

\title{Mild solutions to the 3D-Boussinesq system with weakened initial temperature}
\author{Pedro Gabriel Fern\'andez-Dalgo\footnote{Escuela de Ciencias Físicas y Matemáticas, Universidad de Las Américas, Vía a Nayón, C.P.170504, Quito, Ecuador} \footnote{corresponding author: pedro.fernandez.dalgo@udla.edu.ec} \  and  Oscar Jarr\'in \footnote{Escuela de Ciencias Físicas y Matemáticas, Universidad de Las Américas, Vía a Nayón, C.P.170504, Quito, Ecuador} \footnote{oscar.jarrin@udla.edu.ec}}
\date{}

\maketitle

\begin{abstract}
\noindent
In this research, the Cauchy problem of the 3D viscous Boussinesq system is studied considering an initial temperature with \emph{negative Sobolev regularity}. Precisely,  we construct local in time mild solutions to this system  where the  temperature term belongs to Sobolev spaces of negative order.  Our main contribution is to show how the coupled structure of the Boussinesq system allows us to considerably weaken the regularity in the temperature term.
\end{abstract}
 


\section{Introduction}
We consider the incompressible three dimensional Boussinesq system, which describes the dynamics of a viscous incompressible fluid with heat exchanges \cite{Pedlosky,Salmon}. Mathematicaly, this system  couples  the Navier-Stokes equations and the equations of thermodynamics as follows: 
\begin{equation}\label{Boussinesq}
\begin{cases}
\partial_t \vu - \Delta \vu + (\vu \cdot \vN)\vu +\vN P - \theta \ve_3=0, \\
\partial_t \theta - \Delta \theta + \vu \cdot \vN \theta=0, \\
\text{div}(\vu)=0, \\
\vu(0,\cdot)=\vu_0, \quad \theta(0,\cdot)=\theta_0.
\end{cases}
\end{equation}
Here, $\vu: [0,+\infty)\times \Rt \to \Rt$ is the velocity of the fluid,  $P:[0,+\infty)\times \Rt \to \R$ is the pressure and $\theta:[0,+\infty)\times \Rt \to \R$ is the temperature. Moreover, $\ve_3=(0,0,1)$ is the third orthonormal vector in the canonical basis of $\Rt$, while $\vu_0: \Rt \to \Rt$ and $\theta_0 : \Rt \to \R$ denote the initial (divergence-free) velocity and the initial temperature respectively. With a minor loss of generality, we have set all the physical constants equal to one. 

\medskip

\noindent 
The Boussinesq system (\ref{Boussinesq}) was studied in the $L^p$ setting in \cite{CaDi} and weak solutions were constructed in \cite{YKa,Mo}.  In recent years this model have taken interest in the fluids mechanics community. On the one hand, in view of its proximity to the incompressible axisymmetric Euler system with swirl, the \emph{partially viscous} Boussinesq system in dimension two is of great interest, it was addressed for example in \cite{ KLN, LPZ, HoLi}. On the other hand, in \cite{DaPa} the authors show that the classical results for the standard Navier-Stokes system  remain true for the Boussinesq system even considering a null viscosity for the temperature term:
\begin{equation}\label{Boussinesq-Null-viscosity}
\begin{cases}\vspace{2mm}
\partial_t \vu - \Delta \vu + (\vu \cdot \vN)\vu +\vN P - \theta \ve_3=0, \quad \text{div}(\vu)=0, \\
\partial_t \theta + \vu \cdot \vN \theta=0, 
\end{cases}
\end{equation}
and considering \emph{compatible} regularity conditions for both the initial velocity and the initial temperature. Precisely, in Theorem 1.1 of the paper \cite{DaPa}, global in time weak solutions are constructed under the well-known  initial condition for the  velocity $\vu_0 \in L^2(\Rt)$, and the initial condition for the  temperature $\theta_0 \in L^p(\mathbb R^3)$, $6/5 < p \leq  2$. 

\medskip

\noindent
Concerning the case of initial data belonging to Sobolev spaces, and recalling the embedding $L^p(\Rt)\subset \dot{H}^{-s}(\Rt)$, with $6/5<p\leq 2$ and $0\leq s <1$,  this last condition on the initial temperature suggests to look for solutions of the system (\ref{Boussinesq}) with  $\theta_0 \in \dot{H}^{-s}(\Rt)$.

\medskip


\noindent
In this article, we address the stronger setting of mild solutions and we will show the existence and uniqueness of mild solutions for the three dimensional viscous Boussinesq system (\ref{Boussinesq}) with initial temperature  in Sobolev spaces of \emph{negative order}. In this context, we show that the coupled structure of the  system (\ref{Boussinesq}) allows us to weaken the regularity in the temperature term. 

\medskip

\noindent
Recall that \emph{mild} solutions to the Boussinesq system (\ref{Boussinesq}) solve the following coupled system of integral equations:
\begin{equation}\label{Mild-u}
\begin{split}
\vu(t,\cdot)=&\, e^{t\Delta} \vu_0 - \int_{0}^{t} e^{(t-\tau)\Delta} \P \left( (\vu \cdot \vN)\vu\right)(\tau,\cdot)d \tau \\
&\,+ \int_{0}^{t} e^{(t-\tau)\Delta} \P \left(\theta \ve_3 \right)(\tau,\cdot)d \tau,
\end{split}
\end{equation}
\begin{equation}\label{Mild-theta}
\theta(t,\cdot)=e^{t\Delta} \theta_0 - \int_{0}^{t} e^{(t-\tau)\Delta}\vu \cdot \vN \theta (\tau,\cdot)d \tau, 
\end{equation}
where, for the heat kernel  $h_t$ we denote $e^{t\Delta}\varphi = h_t \ast \varphi$, and $\P$ stands for the Leray's projector. 

\medskip

\noindent
Our key remark is the fact that expression (\ref{Mild-u}) formally verifies the equation
\[ \partial_t \vu - \Delta \vu + (\vu \cdot \vN)\vu +\vN P = \theta \ve_3, \quad \text{div}(\vu)=0,\]
where  $\theta \ve_3$ acts as a \emph{source term} of the classical Navier-Stokes equations. Fujita and  Kato's theory of mild solutions in Sobolev spaces \cite{FuKa},\cite[Chapter $5.2$]{BaCheDan},\cite[Chapter $7.4$]{PLe}  shows that this equation can be studied by considering initial velocities $\vu_0 \in H^r(\Rt)$ with $r\geq 1/2$, and \emph{source terms} belonging to the space $\dot{H}^{r-1}(\Rt)$. In particular, the limit value $r=1/2$ suggests the minimal regularity condition $\theta \in \dot{H}^{-1/2}(\Rt)$.

\medskip

\noindent
{\bf Main results.} Coming back to the Boussinesq system (\ref{Boussinesq}), and motivated by this last remark, for $r\geq 1/2$ we shall assume that $\vu_0 \in H^r(\Rt)$, with $div(\vu_0)=0$, and for $0\leq s \leq 1/2$ we shall assume that $\theta_0 \in \dot{H}^{-s}(\Rt)$. We thus look for conditions on $r$ and $s$ to prove the existence of local in time solutions. 

\medskip

\noindent
In order to simplify the statement of our next theorem, we will say that the system (\ref{Boussinesq}) is locally solved in the space $H^r(\Rt)\times \dot{H}^{-s}(\Rt)$  if for any initial data $(\vu_0, \theta_0)\in H^r(\Rt)\times \dot{H}^{-s}(\Rt)$ (with $div(\vu_0)=0$) there exists a time $T_0>0$, depending on $\vu_0$ and $\theta_0$, and there exists a couple
\begin{equation*}
	\vu\in L^{\infty}([0,T_0], H^r(\Rt))\cap L^{2}([0,T_0],\dot{H}^{r+1}(\Rt)),
	\end{equation*}
		\begin{equation*} 
	\theta \in L^{\infty}([0,T_0], \dot{H}^{-s}(\Rt))\cap L^{2}([0,T_0],\dot{H}^{-s+1}(\Rt)),
	\end{equation*}
	which is a solution to (\ref{Mild-u})-(\ref{Mild-theta}). In this setting, our main result reads as follows:

\begin{Theorem}\label{Th1} Let $1/2\leq r < 2 $ and  let $0\leq s \leq  1/2$.
\begin{enumerate}
    \item If $s<1/2<r$ and 
    	\begin{equation}\label{Relationship-s-r}
	1 \leq  s+r<2,
	\end{equation}
 then  the Boussinesq system (\ref{Mild-u})-(\ref{Mild-theta}) is locally solved  in the space $H^r(\Rt) \times \dot{H}^{-s}(\Rt)$. Moreover, the obtained solution  is the unique one. 
  \medskip
 
 \item In the limit case $s=1/2$, and for $1/2\leq r \leq 1$,  the Boussinesq system (\ref{Mild-u})-(\ref{Mild-theta}) is locally solved   in the space  $H^r(\Rt) \times \dot{H}^{-1/2}(\Rt)$. 
%
%
%
\end{enumerate} 
\end{Theorem}	

\medskip

\noindent
The following comments are in order. Existence of local in time solutions to  (\ref{Mild-u})-(\ref{Mild-theta})  will be obtained by a fixed point argument joint with some sharp regularizing effects of the heat kernel. However,  equation (\ref{Mild-theta})
imposes new  defies.  On the one hand,  we deal  with Sobolev spaces of negative order and, on the other hand, the term $\theta \vu$ is more difficult to treat due the fact $\theta$ and $\vu$ have different regularity properties.  To overcome these difficulties, we have considered two cases of the parameters $r$ and $s$.

\medskip

\noindent
In the first point of Theorem \ref{Th1}, we consider the case $s<1/2<r$ with the additional relationship (\ref{Relationship-s-r}). Observe that this relationship constraints $s$ in function of $r$. Precisely,  by the lower bound in (\ref{Relationship-s-r})  we have 
\begin{equation}\label{Relation-r-s-2}
 -s \leq r-1,    
\end{equation}
and then $-s$ in contained in the interval $(-1/2,r-1)$. This fact shows us how much the regularity  of the initial temperature $\theta_0\in \dot{H}^{-s}(\Rt)$ can be weakened respect to the given regularity  of the initial velocity $\vu_0 \in H^r(\Rt)$. On the other hand, the upper bound in  (\ref{Relationship-s-r}) constraints the parameter $r$ to $r<2$.  We thus observe that the given regularity for initial velocity $\vu_0$ cannot be arbitrary high, when  initial temperature with negative regularity is considered. 

\medskip

\noindent
It is interesting to observe that  the inequality (\ref{Relation-r-s-2})  appears inverted in the recent paper \cite{PFe} on the micropolar system, the reason is the different nature of the coupled part in the equation for $\partial_t \vu$. While the coupled part in the Boussinesq system is $-\theta {\bf e}_3$, the coupled part in the micropolar system is $\nabla \wedge \omega$, where $\omega$ is the coupled variable.

\medskip

\noindent
Getting back to (\ref{Relation-r-s-2}), we observe that when $r$ goes to $1/2$ then $-s$ tends to the limit value $-1/2$. Thus, in the second  point of Theorem \ref{Th1}, we focus on the (more delicate) case of minimal regularity conditions for the initial temperature: $\theta_0 \in \dot{H}^{-1/2}(\Rt)$,  and we show the existence of a local solution to (\ref{Mild-u})-(\ref{Mild-theta}) where the initial velocity $\vu_0$ verifies $\vu_0 \in H^r(\Rt)$ with $1/2\leq r \leq 1$. 

\medskip 

\noindent
In the range $1/2\leq r \leq 1$, the limit points $r=1/2$ and $r=1$ are of particular interest. On the one hand, when $r=1/2$ we deal with the limit case of the relationship (\ref{Relation-r-s-2}). On the other hand, the value $r=1$ seems to be the maximal one for which we can prove the local-well posedness of (\ref{Mild-u})-(\ref{Mild-theta}) in $H^r(\Rt)\times \dot{H}^{-1/2}(\Rt)$. 

\medskip

\noindent
When $s=1/2$, the uniqueness  of solutions seems more complicate to treat. In fact, the methods used in the previous case when $s<1/2$ are not longer valid when studying the  term $\theta \vu$ in equation (\ref{Mild-theta}) with  $\theta \in \dot{H}^{-1/2}(\Rt)$. As pointed out in \cite{BranHe}, this particular term makes more difficult to study the uniqueness issue of solutions to the Boussinesq system, and additional regularity  conditions on the temperature term are required to obtain partial uniqueness  results.

\medskip

\noindent
Following these ideas, we are able to obtain the next (partial) uniqueness result.
\begin{Proposition}\label{Prop-Uniqueness} Withing the setting of the second point in Theorem \ref{Th1}, assume that we have two solutions $(\vu_1,\theta_1)$ and $(\vu_2,\theta_2)$ to the Boussinesq system (\ref{Boussinesq}) associated with the same initial data, and such that 
\begin{equation*}
	\vu_i \in L^{\infty}([0,T_0], H^{1/2}(\Rt))\cap L^{2}([0,T_0],\dot{H}^{3/2}(\Rt)),
	\end{equation*}
		\begin{equation*} 
	\theta_i \in L^{\infty}([0,T_0], \dot{H}^{-1/2}(\Rt))\cap L^{2}([0,T_0],\dot{H}^{1/2}\cap \dot{W}^{1,3}(\Rt)),
	\end{equation*}
 for $i=1,2$. Then $(\vu_1,\theta_1)=(\vu_2,\theta_2)$.
\end{Proposition}

\noindent
As noticed, uniqueness of solutions is ensured under  the additional regularity condition   $\theta \in L^2_t \dot{W}^{1,3}_x$. Moreover, this result also holds for velocities $\vu$ belonging to the space $L^\infty_t H^r_x \cap L^2_t \dot{H}^{r+1}_x$ with $1/2<r\leq 1$, due to the continuous embedding  $L^\infty_t H^r_x \cap L^2_t \dot{H}^{r+1}_x \subset L^\infty_t H^{1/2}_x \cap L^2_t \dot{H}^{3/2}_x$. 

\medskip

\noindent

\medskip

\noindent
{\bf Related works.} Here we make a short discussion on some previous related works. To this end,  some mentioned results  are not stated in their rigorous form, since it requires a considerably set of highly technical  definition and notation. But, we shall emphasize  their main features concerning the \emph{regularity} of the spaces involved. For all the technical details we refer to the articles cited below. 

\medskip

\noindent
Previous studies on the Boussinesq systems (\ref{Boussinesq}) and (\ref{Boussinesq-Null-viscosity}) in functional spaces of negative or null regularity were done in the framework of Besov spaces, generally  defined on  the space $\R^n$ with $n\geq 2$. For our purposes, we shall only focus on the case $n=3$. In \cite{DeCu},  the authors work with (\ref{Boussinesq}) and  essentially consider initial velocities $\vu_0 \in  B^{-1}_{\infty,1}(\Rt)$  and initial temperatures $\theta_0 \in  B^{-1}_{3/2,1}(\Rt)$. Then, it is proven that small data yields local in time solutions to the system (\ref{Boussinesq}), which in essence verify 
\[(\vu,\theta)\in L^\infty_t  (B^{-1}_{\infty,1})_x \cap L^2_t (B^{0}_{\infty,1})_x \times L^\infty_t  (B^{-1}_{3/2,1})_x \cap L^2_t (B^{0}_{3/2,1})_x. \]
The main objective of this work is to perform sharp estimates on the coupling term $\vu \cdot \vN \theta$ to obtain the local well-posedness in low regularity Besov spaces with index $-1$. 

\medskip

\noindent
On the other hand, in \cite[Theorem $1.3$]{DaPa} the authors consider the Boussinesq system (\ref{Boussinesq-Null-viscosity}). They also consider any initial   data $\vu_0, \theta_0\in \dot{B}^{0}_{3,1}(\Rt)$ to construct local in time solutions  
\[(\vu,\theta)\in L^\infty_t  (\dot{B}^{0}_{3,1})_x \cap L^1_t (\dot{B}^{2}_{3,1})_x \times L^\infty_t  (\dot{B}^{0}_{3,1})_x.\]
By scaling properties, we have that  the space $\dot{B}^{0}_{3,1}(\Rt)$  is embedded in $L^3(\Rt)$, which is the well-known scale invariant space for the first equation in (\ref{Boussinesq-Null-viscosity}) involving the velocity $\vu$. 

\medskip

\noindent
Compared with these results, the main difference with our work bases on the fact that we exploit  the coupled structure of the system (\ref{Boussinesq}) (mainly the crossed term $\vu \cdot \vN \theta$) and we use some sharp smoothing effects of the heat kernel (see for instance the third point of  Lemma \ref{heat-estimates-2}  and Lemma \ref{heat-estimates-3} below) to consider   \emph{different} regularity properties for $\vu$ and $\theta$, principally with a considerably  \emph{weakened} regularity on the temperature term $\theta$.  

\medskip

\noindent
{\bf Open questions and future research.} Our first natural question is to look for if Theorem \ref{Th1} still holds for the null viscosity system (\ref{Boussinesq-Null-viscosity}). This is not a trivial fact due to the loss of smoothing effects in equation involving the temperature term $\theta$. 

\medskip

\noindent
On the other hand,  we emphasize that   in contrast to the classical Navier-Stokes equations (when $\theta \equiv 0$), global in time mild  solutions in Sobolev spaces arising from  small initial data (controlled by universal constants) seem  not to be a  trivial issue for both  Boussinesq systems (\ref{Boussinesq}) and (\ref{Boussinesq-Null-viscosity}). In fact, in our case  mild solutions  to (\ref{Boussinesq}) are constructed by an iterative fixed point argument given in Lemma \ref{Picard} below, and  the main difficulty focuses on the linear term in equation (\ref{Mild-u}): in the required estimate $\left\| \int_{0}^{t}e^{(t-\tau)\Delta}\P(\theta \ve_3)(\tau,\cdot)d \tau \right\|_{H^r} \leq C_L \|\theta \|_{\dot{H}^{-s}}$, the continuity constant $C_L>0$ \emph{depends} on the time $T$. However,  due to the first technical constraint  in (\ref{Conditions-Picard}), this constant  must be small enough \emph{independent} of the size of initial data, and this fact blocks to apply classical arguments to construct global in time solutions. 

\medskip

\noindent
Coming back to \cite{DaPa},  Theorem $1.4$ yields that local in time solutions to (\ref{Boussinesq-Null-viscosity}), obtained in Theorem $1.3$ from initial data $\vu_0, \theta_0\in \dot{B}^{0}_{3,1}(\Rt)$, can be extended to global ones under the supplementary hypothesis $\vu_0 \in L^{3,\infty}(\Rt)$, $\theta_0 \in L^1(\Rt)$, and with smallness conditions on the quantity $\| \vu_0 \|_{L^{3,\infty}}+\| \theta_0\|_{L^1}$, which are given by universal constants. In future research, we aim to adapt this method to our framework.  Nevertheless, this does not  seem  to be  trivial, principally when handling with Sobolev spaces of negative regularity for the temperature term $\theta$.  

\medskip

\noindent
Finally, in further research we also aim to understood the optimality of the relationship (\ref{Relationship-s-r}) involving the parameters $r$ and $s$. Precisely, we aim to study some ill-posedness issues in the complementary cases $s+r<1$ or $2\leq s+r$.

\medskip

\noindent
{\bf Notation and organization of the article}  The Fourier transform (in the spatial variable) of a function $f$ is denoted by $\widehat{f}$, while $\mathcal{F}^{-1}_x(f)$ stands for the inverse Fourier transform.  Moreover, in our estimates we shall use a generic constant $C>0$ which may change from one line to another. 

\medskip

\noindent
 In Section \ref{Sec:Preliminaries} we collect all the technical lemmas which we shall use later. Section \ref{Sec:Proof-Main-Th} is devoted to prove Theorem \ref{Th1}, while in Section \ref{Sec:Proof-Prop-Uniq} we give a proof of Proposition \ref{Prop-Uniqueness}.

\section{Preliminaries}\label{Sec:Preliminaries}
In this section, we summarize some known estimates concerning smoothing effects of the heat kernel. For the sake of completeness, in some cases we give a short proof of these statements.
\begin{Lemma}\label{heat-estimates} Let $s_1\in \R$ and $s_2\geq 0$  be two real numbers. There exists a constant $C>0$, which depends on $s_1$ and $s_2$, such that  for any $t>0$ we have:
\[ \| e^{t\Delta} f \|_{H^{s_1+s_2}} \leq C (1+ t^{-\frac{s_2}{2}})\| f \|_{H^{s_1}}. \]
\end{Lemma}
The proof of this lemma is straightforward. It follows from direct computations in the Fourier variable and the well-know identity $\widehat{e^{t\Delta} f}(\xi)=e^{-t|\xi|^2}\widehat{f}(\xi)$. 
\begin{Lemma}\label{heat-estimates-2} Let $f \in L^{2}([0,+\infty),L^2(\Rt)\cap \dot{H}^{s_1}(\Rt))$, with $s_1 \in \R$. Define the function 
	\[ F(t,\cdot)=\int_{0}^{t}e^{(t-\tau)\Delta}f(\tau,\cdot)d \tau. \]
Then, the following estimates hold:
\begin{enumerate}
	\item For all $t>0$ we have $\| \vN F(t,\cdot)\|_{L^2} \leq C \| f \|_{L^{2}_{t}L^{2}_{x}}$.
	\item We have $\| \Delta F \|_{L^{2}_{t}L^{2}_{x}} \leq C \| f \|_{L^{2}_{t}L^{2}_{x}}$.
 \item Let  $s_1 \in \R$ and $1<s_2<2$. Define $p=\frac{2}{s_2-1}$, which verifies $2<p<+\infty$. Then we have the estimate:
 \[ \| F \|_{L^{p}_t \dot{H}^{s_1+s_2}_{x}}\leq C\, \| f \|_{L^{2}_t \dot{H}^{s_1}_{x}}. \]  
\end{enumerate}		
\end{Lemma}	
\pv The first point and the second point are well-known facts, see \cite[Lemma $7.2$]{PLe} for a proof. The third point follows from the previous ones: we write 
\begin{equation*}
\begin{split}
\| F \|_{L^{p}_{t}\dot{H}^{s_1+s_2}_{x}}=&\, \left\| \int_{0}^{t} e^{(t-\tau)\Delta} (-\Delta)^{(s_1+s_2)/2}f(\tau,\cdot) d \tau \right\|_{L^{p}_{t}L^{2}_{x}} \\
=&\, \left\| \int_{0}^{t} e^{(t-\tau)\Delta} (-\Delta)^{\frac{s_1}{2}}f(\tau,\cdot) d \tau \right\|_{L^{p}_{t}\dot{H}^{s_2}_{x}}.
\end{split}
\end{equation*}
For simplicity, denote $G=\int_{0}^{t} e^{(t-\tau)\Delta} g(\tau,\cdot) d \tau$, with $g=(-\Delta)^{\frac{s_1}{2}}f$. We apply  interpolation inequalities in homogeneous Sobolev spaces with  the parameter $\sigma=-s_2+2 \in (0,1)$, hence we have $1-\sigma=s_2 -1$. Then, we use the first and second point stated above to get 
\begin{equation*}
    \| G \|_{L^{p}_{t}\dot{H}^{s_2}_{x}} \leq C\, \| G\|^{\sigma}_{L^{\infty}_{t}\dot{H}^{1}_{x}}\, \| G\|^{1-\sigma}_{L^{2}_{t}\dot{H}^{2}_{x}}\leq C\, \| g \|_{L^{2}_{t}L^{2}_{x}}= C\,\| f \|_{L^{2}_{t}\dot{H}^{s_1}_{x}}. \qquad \blacksquare
\end{equation*}

\noindent
Our last lemma is essentially proven in \cite[Theorem $5.4$]{BaCheDan}. However, we will state it in a more general version adapted to our needs in this article. 
\begin{Lemma}\label{heat-estimates-3} Let $s_1 \in \R$ and let $s_1 < s_2 < s_1 +1$. Define $p=\frac{2}{s_2-s_1}$, which verifies $2<p<+\infty$. 

\medskip

\noindent
For all $\varepsilon>0$ there exists a quantity $R_\varepsilon>0$ such that 
\begin{equation*}
    \left\| e^{t \Delta} f  \right\|_{L^{p}_{t}\dot{H}^{s_2}_{x}}  \leq \frac{\varepsilon}{2}+ (R^2_\varepsilon T)^{1/p} \| f \|_{\dot{H}^{s_1}}.
\end{equation*}
\end{Lemma}
\pv For a parameter $\kappa>0$ (which will be set later) we write
\begin{equation*}
\begin{split}
\left\| e^{t \Delta} f  \right\|_{L^{p}_{t}\dot{H}^{s_2}_{x}} \leq &\, \left\| \mathcal{F}^{-1}_{x} \Big(e^{- t |\xi|^2} \mathds{1}_{|\xi|\geq \kappa}(\xi)\, \widehat{f}\Big)  \right\|_{L^{p}_{t}\dot{H}^{s_2}_{x}}\\
&\, + \left\| \mathcal{F}^{-1}_{x} \Big(e^{- t |\xi|^2} \mathds{1}_{|\xi|< \kappa}(\xi)\, \widehat{f}\Big)  \right\|_{L^{p}_{t}\dot{H}^{s_2}_{x}}. 
\end{split}
\end{equation*}
In order to estimate the first term, recall the the identity $p=\frac{2}{s_2-s_1}$ and the relationship $s_1 < s_2 < s_1+1$. 
Using interpolation inequalities in homogeneous Sobolev spaces with $\sigma=s_1+1-s_2 \in (0,1)$ and $1-\sigma=s_2-s_1$, we have 
\begin{equation*}
\begin{split}
&\, \left\| \mathcal{F}^{-1}_{x} \Big(e^{- t |\xi|^2} \mathds{1}_{|\xi|\geq \kappa}(\xi)\, \widehat{f}\Big)  \right\|_{L^{p}_{t}\dot{H}^{s_2}_{x}}\\
 \leq &\,  C \left( \int_{0}^{T}\left\|  \mathcal{F}^{-1}_{x} \Big(e^{- t |\xi|^2} \mathds{1}_{|\xi|\geq \kappa}(\xi)\, \widehat{f}\Big) \right\|^{\sigma p}_{\dot{H}^{s_1}} \, \left\|  \mathcal{F}^{-1}_{x} \Big(e^{- t |\xi|^2} \mathds{1}_{|\xi|\geq \kappa}(\xi)\, \widehat{f}\Big) \right\|^{(1-\sigma) p}_{\dot{H}^{s_1+2}}  dt\right)^{1/p} \\
 \leq &\, C \left\|  \mathcal{F}^{-1}_{x} \Big(e^{- t |\xi|^2} \mathds{1}_{|\xi|\geq \kappa}(\xi)\, \widehat{f}\Big)\right\|^{\sigma}_{L^{\infty}_{t}\dot{H}^{s_1}_{x}}\, \left\| \mathcal{F}^{-1}_{x} \Big(e^{- t |\xi|^2} \mathds{1}_{|\xi|\geq \kappa}(\xi)\, \widehat{f}\Big) \right\|^{1-\sigma}_{L^{2}_{t}\dot{H}^{s_1+1}_{x}} \\ 
\leq &\, C   \left\|  \mathcal{F}^{-1}_{x} \Big(e^{- t |\xi|^2} \mathds{1}_{|\xi|\geq \kappa}(\xi)\, \widehat{f}\Big)\right\|_{L^{\infty}_{t}\dot{H}^{s_1}_{x}}+ C\, \left\| \mathcal{F}^{-1}_{x} \Big(e^{- t |\xi|^2} \mathds{1}_{|\xi|\geq \kappa}(\xi)\, \widehat{f}\Big) \right\|_{L^{2}_{t}\dot{H}^{s_1+1}_{x}}. 
\end{split}
\end{equation*}
By well-known properties of the heat kernel, each term above is controlled by the quantity $\left\|  \mathcal{F}^{-1}_{x} \Big( \mathds{1}_{|\xi|\geq \kappa}(\xi)\, \widehat{f}\Big)\right\|_{\dot{H}^{s_1}}$,  and we get
\[ \left\| \mathcal{F}^{-1}_{x} \Big(e^{- t |\xi|^2} \mathds{1}_{|\xi|\geq \kappa}(\xi)\, \widehat{f}\Big)  \right\|_{L^{p}_{t}\dot{H}^{s_2}_{x}} \leq C\,  \left\|  \mathcal{F}^{-1}_{x} \Big( \mathds{1}_{|\xi|\geq \kappa}(\xi)\, \widehat{f}\Big)\right\|_{\dot{H}^{s_1}}.\]
Thereafter,  since $f \in \dot{H}^{s_1}(\Rt)$, for $\varepsilon>0$ we can set $\kappa = R_\varepsilon>0$ big enough such that 
\[ C \, \left\| \mathcal{F}^{-1}_{x} \Big(\mathds{1}_{|\xi|\geq R_\varepsilon}(\xi)\, \widehat{f}\Big)\right\|_{\dot{H}^{s_1}} \leq \frac{\varepsilon}{2}.\]

\medskip

\noindent
For the second term, using again the identity $p=\frac{2}{s_2-s_1}$ we write
\begin{equation*}
  \begin{split}
 &\, \left\| \mathcal{F}^{-1}_{x} \Big(e^{- t |\xi|^2} \mathds{1}_{|\xi|< R_\varepsilon}(\xi)\, \widehat{f}\Big)  \right\|_{L^{p}_{t}\dot{H}^{s_2}_{x}}\\
 = &\,  \left( \int_{0}^{T} \|\, |\xi|^{s_2}\,  e^{- t |\xi|^2} \mathds{1}_{|\xi|< R_\varepsilon}(\xi)\, \widehat{f} \|^{p}_{L^2} dt \right)^{1/p} \\
= &\, \left( \int_{0}^{T} \|\, |\xi|^{s_2-s_1}\,  e^{- t |\xi|^2} \mathds{1}_{|\xi|< R_\varepsilon}(\xi)\, |\xi|^{s_1}\widehat{f} \|^{p}_{L^2} dt \right)^{1/p}\\
\leq & \, (R^{p(s_2-s_1)}_{\varepsilon} T)^{1/p} \| f \|_{\dot{H}^{s_1}}=(R^{2}_{\varepsilon} T)^{1/p} \| f \|_{\dot{H}^{s_1}}. \qquad \qquad \blacksquare
  \end{split}  
\end{equation*}
\section{Proof of Theorem \ref{Th1}}\label{Sec:Proof-Main-Th}
In equations (\ref{Mild-u}) and (\ref{Mild-theta}),  we have  a bilinear term 
\begin{equation*}
    B\Big( (\vu, \theta) \, , \,(\tilde \vu, \tilde \theta) \Big) = \Big( B_1( \, (\vu, \theta) \, , \,(\tilde \vu, \tilde \theta)\,), \, B_2 (\,  (\vu, \theta) \, , \,(\tilde \vu, \tilde \theta) \,) \Big),
\end{equation*}
where $B_1$ involves only $\vu$ and $\tilde \vu$, while $B_2$ involves $\vu$ and $\tilde \theta$, as follows:
\begin{align*}
    B_1( \, (\vu, \theta) \, , \,(\tilde \vu, \tilde \theta) \, ) &= - \int_{0}^{t} e^{(t-\tau)\Delta} \P \left( (\vu \cdot \vN) \tilde \vu\right)(\tau,\cdot)d \tau , \\
    B_2( \, (\vu, \theta) \, , \,(\tilde \vu, \tilde \theta) \, ) &= - \int_{0}^{t} e^{(t-\tau)\Delta}\vu \cdot \vN \tilde \theta (\tau,\cdot)d \tau.
\end{align*}
\noindent
Moreover, we have a liner term involving $\theta$,
\begin{equation*}
    L( \, (\vu, \theta) \,) = ( \, \, L_1( \, (\vu, \theta)  \, ), \, L_2( \, (\vu, \theta)  \,) \, \, ),
\end{equation*}
with
\begin{equation*}
    L_1( \, (\vu, \theta)  \, ) = \int_{0}^{t} e^{(t-\tau)\Delta} \P \left(\theta \ve_3 \right)(\tau,\cdot)d \tau \quad \text{and} \quad   L_2 \equiv 0.
\end{equation*}

\noindent
Thus, letting $e=(\vu,\theta)$ and $e_0 = (\, e^{t\Delta} \vu_0, \, e^{t\Delta} \theta_0 \, )$, the whole system (\ref{Mild-u})-(\ref{Mild-theta}) is written of the form
\begin{equation}\label{Fixed-point}
e=e_0+B(e,e)+L(e),
\end{equation}
and to construct a solution we use the following version of the Picard's iteration scheme:
\begin{Lemma}\label{Picard} Let $(E, \| \cdot \|_E)$ be a Banach space and let $e_0 \in E$ be an initial datum. We set $\| e_0 \|_E \leq \delta$. Moreover, let $B: E \times E \to E$ be a bilinear form and let $L:E \to E$ be a linear form, which, for all $e,f \in E$ verify
\begin{equation*}
\| B(e,f) \|_{E} \leq C_B \| e \|_E\, \| f \|_E \ \ \mbox{and} \ \ \| L(e)\|_{E} \leq C_L \| e \|_E. 
\end{equation*}
If the constants $C_B>0$ and $C_L>0$ satisfy the relationships:
\begin{equation}\label{Conditions-Picard}
0<C_L<\frac{1}{3}, \ \ 0< 9 C_B \, \delta <1 \ \ \mbox{and} \ \ C_L + 6 C_B\, \delta <1,
\end{equation}
then equation \ref{Fixed-point} has a solution $e$, which is  uniquely defined by $\| e \|_E \leq 3\delta$. 
\end{Lemma}	
For a proof, we refer to \cite{ChYa} (proof oh Theorem 3.2 in Appendix). Observe that the last inequality in \eqref{Conditions-Picard} is consequence of the two first inequalities.

\medskip

\noindent
Now, to prove Theorem \ref{Th1}  we shall consider the cases $s<1/2<r$ and $s=1/2, \ 1/2 \leq r \leq 1$ separately. 
\subsection{The case $0\leq s<1/2<r<2$.}
Let $T>0$ a time. For the sake of simplicity we shall denote the Banach spaces
\[ E_1= L^{\infty}([0,T], H^r(\Rt))\cap L^{2}([0,T],\dot{H}^{r+1}(\Rt)), \]
and 
\[ E_2= L^{\infty}([0,T], \dot{H}^{-s}(\Rt))\cap L^{2}([0,T],\dot{H}^{-s+1}(\Rt)),\]
with the usual norms 
\[ \| f\|_{E_1} = \sup_{0\leq t \leq T} \| f(t,\cdot)\|_{H^r}+ \left(\int_{0}^{T}\| f(t,\cdot)\|^{2}_{\dot{H}^{r+1}} dt\right)^{\frac{1}{2}}, \]
and
\[ \| f\|_{E_2} = \sup_{0\leq t \leq T} \| f(t,\cdot)\|_{\dot{H}^{-s}}+ \left(\int_{0}^{T}\| f(t,\cdot)\|^{2}_{\dot{H}^{-s+1}} dt\right)^{\frac{1}{2}}, \]
respectively. Here the homogeneous Sobolev spaces are defined as the closure of the test functions with respect to the homogeneous Sobolev seminorm. We will use the Picard's iteration scheme (given in Lemma \ref{Picard}) to construct a local in time solution $(\vu,\theta)\in E_1 \times E_2$ of the coupled  system (\ref{Mild-u})-(\ref{Mild-theta}).

\medskip

\noindent
In these equations, terms involving initial data are simple to estimate and we have 
\begin{equation}\label{Initial-data}
\| e^{t\Delta} \vu_0 \|_{E_1} \leq C \| \vu_0 \|_{H^r}, \quad \| e^{t\Delta} \theta_0 \|_{E_2} \leq C \| \theta_0 \|_{\dot{H}^{-s}}.
\end{equation}
Moreover, for the bilinear term in equation (\ref{Mild-u}) we have the well-known estimate
\begin{equation}\label{Bilinear-N-S}
\left\| \int_{0}^{t} e^{(t-\tau)\Delta} \P \left( (\vu \cdot \vN)\vv\right)(\tau,\cdot)d \tau \right\|_{E_1} \leq C T^{\frac{1}{4}\min(1, 2r-1)}\,\| \vu \|_{E_1} \| \vv \|_{E_1}, 
\end{equation}
see for instance \cite[Theorem $7.3$]{PLe}. 

\medskip

\noindent
Thus, the novelty in this proof is to use the information $\vu \in E_1$ and $\theta \in E_2$ to perform sharp estimates on the term $\ds{\int_{0}^{t} e^{(t-\tau)\Delta} \P \left(\theta \ve_3 \right)(\tau,\cdot)d \tau}$ in equation (\ref{Mild-u}), and the term $\ds{\int_{0}^{t} e^{(t-\tau)\Delta} \vu \cdot \vN \theta (\tau,\cdot)d \tau}$ in equation (\ref{Mild-theta}).  For the sake of clearness, we  state these estimates in the following set of technical lemmas. 
\begin{Lemma} We have 
\begin{equation}\label{Linear1}
\left\|  \int_{0}^{t} e^{(t-\tau)\Delta} \P \left(\theta \ve_3 \right)(\tau,\cdot)d \tau \right\|_{E_1} \leq C \left( T+T^{\frac{2-(r+s)}{2}}\right) \| \theta \|_{E_2},
\end{equation}
where   the upper bound in (\ref{Relationship-s-r})  yields that $2-(r+s)>0$. 
\end{Lemma}
\pv  To control the first term in the norm $\| \cdot \|_{E_1}$, for $0<t\leq T$ fixed we write
\begin{equation*}
\begin{split}
\left\| \int_{0}^{t} e^{(t-\tau)\Delta} \P \left(\theta \ve_3 \right)(\tau,\cdot)d \tau \right\|_{H^r}  \leq \int_{0}^{t} \left\| e^{(t-\tau)\Delta} (\theta \ve_3)(\tau,\cdot) \right\|_{H^r} d \tau.  
\end{split} 
\end{equation*} 
Then, we apply  Lemma \ref{heat-estimates} (with $s_1=-s$ and $s_2=r+s$) to write
\begin{equation*}
    \begin{split}
\int_{0}^{t} \left\| e^{(t-\tau)\Delta} (\theta \ve_3)(\tau,\cdot) \right\|_{H^r} d \tau \leq & \, C\int_{0}^{t} (1+ (t-\tau)^{-\frac{r+s}{2}})\| \theta(\tau,\cdot)\|_{H^{-s}} d \tau\\
\leq & \, C\int_{0}^{t} (1+ (t-\tau)^{-\frac{r+s}{2}})\| \theta(\tau,\cdot)\|_{\dot{H}^{-s}} d \tau \\
\leq & \, C (t+t^{1-\frac{r+s}{2}}) \| \theta \|_{L^{\infty}_{t}\dot{H}^{-s}_{x}},
\end{split}
\end{equation*}
where we have used the fact that $ s > 0$.
We then obtain
\begin{equation}\label{Linear01}
   \sup_{0\leq t \leq T}  \left\| \int_{0}^{t} e^{(t-\tau)\Delta} \P \left(\theta \ve_3 \right)(\tau,\cdot)d \tau \right\|_{H^r} \leq C(T+T^{\frac{2-(r+s)}{2}}) \| \theta \|_{E_2}.
\end{equation}
To control the second term in the norm $\| \cdot \|_{E_1}$, we write 
\begin{equation*}
\begin{split}
&\left\| \int_{0}^{t}e^{(t-\tau)\Delta}\P (\theta \ve_3)(\tau,\cdot)d \tau \right\|_{L^{2}_{t}\dot{H}^{r+1}_{x}}\\
=& \left\| (-\Delta)^{\frac{r+1}{2}} \left( \int_{0}^{t}e^{(t-\tau)\Delta}\P(\theta \ve_3)(\tau,\cdot)d \tau \right)  \right\|_{L^{2}_{t}L^{2}_{x}} \\
=& \left\| \Delta  \left( \int_{0}^{t}e^{(t-\tau)\Delta} \P (-\Delta)^{\frac{r-1}{2}} (\theta \ve_3)(\tau,\cdot)d \tau \right)  \right\|_{L^{2}_{t}L^{2}_{x}}, 
\end{split}
\end{equation*}
and  by the second point of Lemma \ref{heat-estimates-2}  we have 
\begin{equation*}
\left\| \Delta  \left( \int_{0}^{t}e^{(t-\tau)\Delta} \P (-\Delta)^{\frac{r-1}{2}} (\theta \ve_3)(\tau,\cdot)d \tau \right)  \right\|_{L^{2}_{t}L^{2}_{x}} \leq C \| \theta \|_{L^{2}_{t}\dot{H}^{r-1}_{x}}.
\end{equation*}
Then, we use  interpolation inequalities in homogeneous Sobolev spaces by writing $r-1=\sigma(-s) + (1-\sigma)(-s+1)$, hence  $\sigma=2-(r+s)$. 
\begin{Remark}
Note that $0<\sigma \leq 1$ as long as (\ref{Relationship-s-r}) holds.    
\end{Remark}
 We thus obtain
\begin{equation}\label{Estim01}
   \begin{split}
  &  \,  C \left(\int_{0}^{T} \| \theta(\tau,\cdot)\|^{2}_{\dot{H}^{r-1}} d \tau \right)^{\frac{1}{2}}\\
   \leq &\, C \left( \int_{0}^{T} \| \theta(\tau,\cdot)\|^{2\sigma}_{\dot{H}^{-s}}  \| \theta(\tau,\cdot)\|^{2(1-\sigma)}_{\dot{H}^{1-s}} d \tau\right)^{\frac{1}{2}}\\
   \leq &\, C \| \theta \|^{\sigma}_{L^{\infty}_{t}\dot{H}^{-s}_{x}} \left( \int_{0}^{T} \| \theta(\tau,\cdot)\|^{2(1-\sigma)}_{\dot{H}^{1-s}} d \tau\right)^{\frac{1}{2}}.
   \end{split}
\end{equation}
In order to control the last integral,  we apply H\"older inequalities with $p=\frac{1}{1-\sigma}$ and $q=\frac{1}{\sigma}$ (hence we have $1=1/p+1/q$) to obtain
\begin{equation}\label{Estim02}
\begin{split}
\left( \int_{0}^{T} \| \theta(\tau,\cdot)\|^{2(1-\sigma)}_{\dot{H}^{1-s}} d \tau\right)^{\frac{1}{2}} \leq &\,  \left( \int_{0}^{T} \| \theta(\tau,\cdot)\|^{2}_{\dot{H}^{1-s}} d \tau\right)^{\frac{1-\sigma}{2}} T^{\frac{\sigma}{2}}\\
=&\,\| \theta \|^{1-\sigma}_{L^{2}_{t}\dot{H}^{1-s}_{x}} \, T^{\frac{\sigma}{2}}.
\end{split}
\end{equation}
Gathering (\ref{Estim01}) and (\ref{Estim02}), and applying the discrete Young inequalities with $p$ and $q$ given above, we  have
\begin{equation}\label{Linear02}
 \left\| \int_{0}^{t}e^{(t-\tau)\Delta}\P(\theta \ve_3)(\tau,\cdot)d \tau \right\|_{L^{2}_{t}\dot{H}^{r+1}_{x}} \leq C T^{\frac{\sigma}{2}}\| \theta \|_{E_2}.   
\end{equation}
Inequality (\ref{Linear1}) directly follows from (\ref{Linear01}) and (\ref{Linear02}), and the identity $\sigma=2-(r+s)>0$.  \finpv

\begin{Lemma} We have
\begin{equation}\label{Bilinear}
\left\| \int_{0}^{t} e^{(t-\tau)\Delta}\vu \cdot \vN \theta (\tau,\cdot)d \tau \right\|_{E_2} \leq C T^{-\frac{s}{4}+\frac{1}{8}} \| \vu \|_{E_1} \| \theta \|_{E_2}, \quad  0<-\frac{s}{4}+\frac{1}{8}. 
\end{equation}
Observe that $0<-\frac{s}{4}+\frac{1}{8}$ as long as $s<\frac{1}{2}$.
\end{Lemma}
\pv  For the first term in $\| \cdot \|_{E_2}$, we have the estimate
\begin{equation}\label{Bilinear1}
\left\| \int_{0}^{t} e^{(t-\tau)\Delta}\vu \cdot \vN \theta (\tau,\cdot)d \tau \right\|_{L^{\infty}_{t}\dot{H}^{-s}_{x}} \leq C T^{-\frac{s}{4}+\frac{1}{8}}\|\vu \|_{E_1}\|\theta\|_{E_2}.   
\end{equation}
Indeed, using the first point in Lemma \ref{heat-estimates-2} we write 
\begin{equation*}
\begin{split}
 &\, \left\| \int_{0}^{t} e^{(t-\tau)\Delta}\vu \cdot \vN \theta (\tau,\cdot)d \tau \right\|_{\dot{H}^{-s}} \\
 =&\,  \left\| (-\Delta)^{-\frac{s}{2}}\left( \int_{0}^{t} e^{(t-\tau)\Delta}\vu \cdot \vN \theta (\tau,\cdot)d \tau  \right)  \right\|_{L^2} \\
 \leq & \, C \left\| \vN\otimes  \left( \int_{0}^{t} e^{(t-\tau)\Delta} (-\Delta)^{-\frac{s}{2}} (\theta \vu)(\tau,\cdot)d \tau  \right)  \right\|_{L^2} \\
 \leq & \, C \| (-\Delta)^{-\frac{s}{2}} (\theta \vu) \|_{L^{2}_{t}L^{2}_{x}} = C \| \theta \vu \|_{L^{2}_{t}\dot{H}^{-s}_{x}}.
 \end{split}
\end{equation*}
To control the last term, by the Product laws in homogeneous Sobolev spaces we have 
\begin{equation}\label{Estim02-tech}
C \| \theta \vu \|_{L^{2}_{t}\dot{H}^{-s}_{x}} \leq C \left( \int_{0}^{T} \| \theta(\tau,\cdot)\|^{2}_{\dot{H}^{-\frac{s}{2}+\frac{3}{4}}}\,\| \vu(\tau,\cdot)\|^{2}_{\dot{H}^{-\frac{s}{2}+\frac{3}{4}}}  d \tau \right)^{\frac{1}{2}}.
\end{equation}
Remark that $0<-\frac{s}{2}+\frac{3}{4} < \frac{3}{2}$ as long as $-\frac{3}{2}<s<\frac{3}{2}$, which is verified for $0\leq s < \frac{1}{2}$.

\medskip

\noindent
In the first term on the right of (\ref{Estim02-tech}),  we  apply interpolation estimates with the relationships  $-\frac{s}{2}+\frac{3}{4}=-s\sigma_1 +(1-s)(1-\sigma_1)$ with $\sigma_1=-\frac{s}{2}+\frac{1}{4}$. Remark that the required control $0<\sigma_1 <1$ is ensured by our assumption $0\leq s < \frac{1}{2}$. Therefore we can write
\[ \| \theta(\tau,\cdot)\|^{2}_{\dot{H}^{-\frac{s}{2}+\frac{3}{4}}} \leq C\, \| \theta(\tau,\cdot)\|^{2\sigma_1}_{\dot{H}^{-s}}\|\theta (\tau,\cdot)\|^{2(1-\sigma_1)}_{\dot{H}^{-s+1}}. \]
Similarly, in the second on the right of (\ref{Estim02-tech}), we use  $-\frac{s}{2}+\frac{3}{4}= \sigma_2 0 + (1-\sigma_2)r$ with $\sigma_2=\frac{s}{2r}-\frac{3}{4r}+1$, where we also have $0<\sigma_2<1$. Indeed, on the one hand, observe that $0<\sigma_2$ as long as $\frac{3}{4}<\frac{s}{2}+r$. However, this last inequality is verified thanks to (\ref{Relationship-s-r}) and the fact that $s<\frac{1}{2}<r$: by the lower bound in (\ref{Relationship-s-r}) we write $\frac{3}{4}\leq \frac{3}{4}s+\frac{3}{4}r$. The inequality $\frac{3}{4}s+\frac{3}{4}r < \frac{s}{2}+r$ is equivalent to the inequality $s < r$, which is ensured by $s<\frac{1}{2}<r$.  On the other hand, observe that $\sigma_2 <1$ as long as $s<\frac{3}{2}$, which is ultimately verified by the assumption  $s<\frac{1}{2}$. Therefore, we write
\[ \| \vu(\tau,\cdot)\|^{2}_{\dot{H}^{-\frac{s}{2}+\frac{3}{4}}} \leq C\, \| \vu(\tau,\cdot)\|^{2\sigma_2}_{L^2}\|\vu (\tau,\cdot)\|^{2(1-\sigma_2)}_{\dot{H}^{r}}. \]
We thus obtain
\begin{equation*}
  \begin{split}
 &\,  C \left( \int_{0}^{T} \| \theta(\tau,\cdot)\|^{2}_{\dot{H}^{-\frac{s}{2}+\frac{3}{4}}}\,\| \vu(\tau,\cdot)\|^{2}_{\dot{H}^{-\frac{s}{2}+\frac{3}{4}}}  d \tau \right)^{\frac{1}{2}}\\
 \leq & \, C \left( \int_{0}^{T} \| \theta(\tau,\cdot)\|^{2\sigma_1}_{\dot{H}^{-s}}\|\theta (\tau,\cdot)\|^{2(1-\sigma_1)}_{\dot{H}^{-s+1}}\, \| \vu(\tau,\cdot)\|^{2\sigma_2}_{L^2}\|\vu (\tau,\cdot)\|^{2(1-\sigma_2)}_{\dot{H}^{r}}d \tau\right)^{\frac{1}{2}}\\
 \leq & \, C \| \theta \|^{\sigma_1}_{L^{\infty}_{t}\dot{H}^{-s}_{x}}\, \| \vu\|^{\sigma_2}_{L^{\infty}_{t}L^{2}_x}\|\vu \|^{1-\sigma_2}_{L^{\infty}_{t}\dot{H}^{r}_{x}}\, \left( \int_{0}^{T} \|\theta (\tau,\cdot)\|^{2(1-\sigma_1)}_{\dot{H}^{-s+1}} d \tau   \right)^{\frac{1}{2}} \\
 \leq &\, C  \|\vu \|_{L^{\infty}_{t}H^{r}_{x}} \, \| \theta \|^{\sigma_1}_{L^{\infty}_{t}\dot{H}^{-s}_{x}}\,\left( \int_{0}^{T} \|\theta (\tau,\cdot)\|^{2(1-\sigma_1)}_{\dot{H}^{-s+1}} d \tau   \right)^{\frac{1}{2}}.
  \end{split}
\end{equation*}
We still need to estimate the last integral. For this we use H\"older estimates (in the variable of time) with $\frac{1}{p}=1-\sigma_1$ and $\frac{1}{q}=\sigma_1$. Moreover, recalling that $0<\sigma_1=-\frac{s}{2}+\frac{1}{4}$  we obtain 
\begin{equation*}
 \begin{split}
&\,  C  \|\vu \|_{L^{\infty}_{t}H^{r}_{x}} \, \| \theta \|^{\sigma_1}_{L^{\infty}_{t}\dot{H}^{-s}_{x}}\,\left( \int_{0}^{T} \|\theta (\tau,\cdot)\|^{2(1-\sigma_1)}_{\dot{H}^{-s+1}} d \tau   \right)^{\frac{1}{2}} \\
\leq &\, C \|\vu \|_{L^{\infty}_{t}H^{r}_{x}} \, \| \theta \|^{\sigma_1}_{L^{\infty}_{t}\dot{H}^{-s}_{x}}\,\left( \int_{0}^{T} \|\theta (\tau,\cdot)\|^{2}_{\dot{H}^{-s+1}} d \tau   \right)^{\frac{1-\sigma_1}{2}} T^{\frac{\sigma_1}{2}}\\
\leq & \, C T^{-\frac{s}{4}+\frac{1}{8}} \| \vu \|_{E_1}\| \theta \|_{E_2},
 \end{split}
\end{equation*}
which yields (\ref{Bilinear1}).

\medskip

\noindent
For the second term in $\| \cdot \|_{E_2}$ we have the estimate 
\begin{equation}\label{Bilinear2}
\left\| \int_{0}^{t} e^{(t-\tau)\Delta}\vu \cdot \vN \theta (\tau,\cdot)d \tau \right\|_{L^{2}_{t}\dot{H}^{-s+1}_{x}} \leq C T^{-\frac{s}{4}+\frac{1}{8}}\|\vu \|_{E_1}\|\theta\|_{E_2}.  
\end{equation}
Indeed, we write 
\begin{equation*}
\begin{split}
&\, \left\| \int_{0}^{t} e^{(t-\tau)\Delta}\vu \cdot \vN \theta (\tau,\cdot)d \tau \right\|_{L^{2}_{t}\dot{H}^{-s+1}_{x}} \\
=&\, \left\| \int_{0}^{t} e^{(t-\tau)\Delta} \text{div}(\theta \vu) (\tau,\cdot)d \tau \right\|_{L^{2}_{t}\dot{H}^{-s+1}_{x}} \\
\leq & \, C  \left\| \int_{0}^{t} e^{(t-\tau)\Delta} (\theta \vu) (\tau,\cdot)d \tau \right\|_{L^{2}_{t}\dot{H}^{-s+2}_{x}}.
\end{split}
\end{equation*}
We use here the second point in Lemma \ref{heat-estimates-2} to write
\begin{equation*}
\begin{split}
&\, C  \left\| \int_{0}^{t} e^{(t-\tau)\Delta} (\theta \vu) (\tau,\cdot)d \tau \right\|_{L^{2}_{t}\dot{H}^{-s+2}_{x}} \\
\leq &\, C \left\| \Delta \left( \int_{0}^{t} e^{(t-\tau)\Delta} (-\Delta)^{-\frac{s}{2}}(  \theta \vu) (\tau,\cdot)d \tau \right) \right\|_{L^{2}_{t}L^{2}_{x}}\\
\leq &\,  C \| (-\Delta)^{-\frac{s}{2}} (\theta \vu) \|_{L^{2}_{t}L^{2}_{x}} = C \| \theta \vu \|_{L^{2}_{t}H^{-s}_{x}},
\end{split} 
\end{equation*}
where the last term  was already estimated above. Estimate (\ref{Bilinear}) now directly follows from estimates (\ref{Bilinear1}) and (\ref{Bilinear2}). \finpv 

\noindent
With  estimates (\ref{Initial-data}), (\ref{Bilinear-N-S}), (\ref{Linear1}) and (\ref{Bilinear})  at our disposal, we set a time $T_0$ small enough in order to satisfy all the set of conditions (\ref{Conditions-Picard}) in Lemma \ref{Picard}, and we obtain a solution $(\vu, \theta)\in E_1 \times E_2$ to the system (\ref{Mild-u})-(\ref{Mild-theta}).

\medskip

\noindent
In addition, this solution is the unique one in the space $E_1 \times E_2$. The proof of this fact is rather standard, so we shall only detail the main ideas. Let $(\vu_1, \theta_1),(\vu_2, \theta_2)$ be two solutions of (\ref{Mild-u})-(\ref{Mild-theta}) arising from the same initial data. We define $\vv=\vu_1 - \vu_2$ and $\eta=\theta_1 - \theta_2$. Then, $(\vv,\eta)$ solves the coupled system:
\begin{equation*}
\begin{split}
\vv(t,\cdot)=&\, - \int_{0}^{t} e^{(t-\tau)\Delta} \P \Big( (\vu_1 \cdot \vN)\vv + (\vv \cdot \vN)\vu_2 \Big)(\tau,\cdot)d \tau \\
&\,+ \int_{0}^{t} e^{(t-\tau)\Delta} \P \left(\eta \ve_3 \right)(\tau,\cdot)d \tau,
\end{split}
\end{equation*}
\begin{equation*}
\eta(t,\cdot)= - \int_{0}^{t} e^{(t-\tau)\Delta} \Big( \vv \cdot \vN \theta_1 + \vu_2 \cdot \vN \eta\Big)  (\tau,\cdot)d \tau.
\end{equation*}

\noindent
On the other hand, we denote by $0\leq T_* \leq T_0$ the maximal time such that $(\vv,\eta)=(0,0)$ on $[0,T_*]\times \Rt$. We shall prove that $T_*=T_0$. For this, we shall assume that $T_*<T_0$ to obtain a contradiction. Indeed, if $T_*<T_0$ let a time $T_*<T_1<T_0$. Then,  we consider  the spaces $E_1$ and $E_2$ on the interval of time $[T_*, T_1]$. By performing again the estimates (\ref{Bilinear-N-S}), (\ref{Linear1} and (\ref{Bilinear}), we obtain
\begin{equation*}
\begin{split}
\| \vv \|_{E_1} \leq  &\, C (T_1-T_*)^{\frac{1}{4}\min(1,2r-1)} \left( \|\vu_1\|_{E_1}+\|\vu_2\|_{E_1}\right) \| \vv \|_{E_1} \\
&\, + C\, \Big((T_1-T_*)+(T_1-T_*)^{\frac{2-(r+s)}{2}}\Big)\| \eta \|_{E_2},
\end{split}
\end{equation*}
and
\begin{equation*}
\| \eta \|_{E_2} \leq C (T_1-T_*)^{-\frac{s}{4}+\frac{1}{8}} \left( \| \vv \|_{E_1}\| \theta_1\|_{E_2}+\|\vu_2 \|_{E_1}\|\eta\|_{E_2}\right).    
\end{equation*}
Here, we can set the time $T_1$ close enough to the time $T_*$ such that $\|\vv \|_{E_1}+\|\eta\|_{E_2}=0$. We obtain that $(\vv,\eta)=(0,0)$ on $[0,T_1)\times \Rt$, which is a contradiction with the definition of the time $T_*$. We thus have $T_*=T_0$.

\subsection{The case $s=1/2$ and $1/2\leq r \leq 1$.}
Due to technical difficulties, and for the sake of clearness, we will divide this case in the following subcases.
\subsubsection{When $s=1/2$ and $r=1/2$.}\label{Sec:s-r-1/2} 
We set
\[ E_1= L^{\infty}([0,T], H^{ \frac{1}{2} }(\Rt))\cap L^{2}([0,T],\dot{H}^{ \frac{3}{2} }(\Rt)), \]
and 
\[ E_2= L^{\infty}([0,T], \dot{H}^{ -\frac{1}{2} }(\Rt))\cap L^{2}([0,T],\dot{H}^{ \frac{1}{2} }(\Rt)),\]
considering the norms
\[ \| f\|_{E_1} = \sup_{0\leq t \leq T} \| f(t,\cdot)\|_{H^{ \frac{1}{2}} }+ \left(\int_{0}^{T}\| f(t,\cdot)\|^{2}_{\dot{H}^{ \frac{3}{2} }} dt\right)^{\frac{1}{2}}, \]
and
\[ \| f\|_{E_2} = \sup_{0\leq t \leq T} \| f(t,\cdot)\|_{\dot{H}^{ -\frac{1}{2} }}+ \left(\int_{0}^{T}\| f(t,\cdot)\|^{2}_{\dot{H}^{ \frac{1}{2} }} dt\right)^{\frac{1}{2}}.\]

\noindent
For the initial data part we have 
\begin{equation}\label{Initial-data2}
\| e^{t\Delta} \vu_0 \|_{E_1} \leq C \| \vu_0 \|_{H^\frac{1}{2}}, \quad \| e^{t\Delta} \theta_0 \|_{E_2} \leq C \| \theta_0 \|_{\dot{H}^{-1/2}}.
\end{equation}
However, in contrast to the previous case $s<\frac{1}{2}<r$,  estimates on the bilinear terms  in equations (\ref{Mild-u}) and (\ref{Mild-theta}) are not longer dependent on the time $T$, and this fact constraints to  consider small initial data $(\vu_0, \theta_0)\in \dot{H}^{1/2}(\Rt) \times \dot{H}^{-1/2}(\Rt)$ in order to construct a local in time solution $(\vu, \theta)\in E_1 \times E_2$ to the system (\ref{Mild-u})-(\ref{Mild-theta}). 

\medskip

\noindent
To overcome this problem (in order to work  with large data) we shall consider different functional spaces to apply the point-fixed argument. Remark that by interpolation inequalities we have the continuous embeddings 
\begin{equation}\label{Embeddings}
E_1 \subset L^4([0,T], \dot{H}^1(\Rt)), \ \ \mbox{and} \ \  E_2 \subset L^4([0,T], L^2(\Rt)).
\end{equation}

\noindent
By (\ref{Initial-data2}) and  (\ref{Embeddings}) we obtain that $e^{t\Delta} \vu_0 \in L^4([0,T], \dot{H}^1(\Rt))$ and $e^{t\Delta} \theta_0 \in L^4([0,T], L^2(\Rt)$, and we thus have 
the following bound, uniformly in $T$,
\begin{equation*}
\begin{split}
 &\| e^{t\Delta} \vu_0 \|_{L^4([0,T], \dot{H}^1(\Rt))} \leq C \| \vu_0 \|_{H^\frac{1}{2}},\\
 & \| e^{t\Delta} \theta_0 \|_{L^4([0,T], L^2(\Rt))} \leq C \| \theta_0 \|_{\dot{H}^{-1/2}}.
\end{split}
\end{equation*}

\noindent
Moreover, in Lemma \ref{heat-estimates-3}  we set  the parameters $(s_1,s_2)=(1/2,1)$   and  $(s_1,s_2)=(-1/2,0)$, hence in both case  we get $p=4$. Then, one have the controls:
\begin{equation}\label{Control-u0}
\left\| e^{t\Delta} \vu_0 \right\|_{L^{4}_{t}\dot{H}^{1}_{x}} \leq \frac{\varepsilon}{2} + C\, (R^{2}_{\varepsilon} \,  T)^{\frac{1}{4}} \| \vu_0 \|_{\dot{H}^{1/2}},   
\end{equation}   
and
\begin{equation}\label{Control-theta0}
 \left\| e^{t\Delta} \theta_0 \right\|_{L^{4}_{t}L^{2}_{x}} \leq \frac{\varepsilon}{2} + C\, (R^{2}_{\varepsilon} \,  T)^{\frac{1}{4}} \| \theta_0 \|_{\dot{H}^{-1/2}}.   
\end{equation}

\medskip

\noindent
Consequently, estimates (\ref{Control-u0}) and (\ref{Control-theta0}) allow us to consider any initial data 
\[ (\vu_0,\theta_0)\in \dot{H}^{1/2}(\Rt) \times \dot{H}^{-1/2}(\Rt),\]
and to still  use the Picard's iteration scheme to construct a (local in time) solution $(\vu, \theta) \in L^4([0,T], \dot{H}^1(\Rt)) \times L^4([0,T], L^2(\Rt) )$ to the coupled system (\ref{Mild-u})-(\ref{Mild-theta}).

\medskip

\noindent As before, we shall  estimate all the terms in (\ref{Mild-u})-(\ref{Mild-theta}). For the bilinear term in equation (\ref{Mild-u}), by \cite[Corollary $5.11$]{BaCheDan} we have the well-known estimate
 \begin{equation}\label{Bilinear-N-S2}
\left\| \int_{0}^{t} e^{(t-\tau)\Delta} \P \left( (\vu \cdot \vN)\vv\right)(\tau,\cdot)d \tau \right\|_{L^{4}_{t}\dot{H}^{1}_{x}} \leq C \| \vu \|_{L^{4}_{t}\dot{H}^{1}_{x}} \| \vv \|_{L^{4}_{t}\dot{H}^{1}_{x}},
\end{equation}   
where the constant $C>0$ does not depend on the time $T$. So, we must focus on the rest of the terms in equations (\ref{Mild-u})-(\ref{Mild-theta}). 

\begin{Lemma} We have
 \begin{equation}\label{Linear1-limit-case}
\left\|  \int_{0}^{t} e^{(t-\tau)\Delta} \P \left(\theta \ve_3 \right)(\tau,\cdot)d \tau \right\|_{L^{4}_{t}\dot{H}^{1}_{x}} \leq C T^{\frac{1}{2}} \| \theta \|_{L^{4}_{t}L^{2}_{x}}.
\end{equation} 
\end{Lemma}
\pv  By the first point in Lemma \ref{heat-estimates-2} we write
\begin{equation*}
\begin{split}
&\left( \int_{0}^{T} \left\| \int_{0}^{t} e^{(t-\tau)\Delta} \P \left(\theta \ve_3 \right)(\tau,\cdot)d \tau \right\|^{4}_{\dot{H}^1}  dt\right)^{\frac{1}{4}}
\leq   \left( \int_{0}^{T} \left\| \theta \right\|^{4}_{L^{2}_{t}L^{2}_{x}}  dt\right)^{\frac{1}{4}} \\
\leq & \, C\, T^{\frac{1}{4}}  \left\| \theta \right\|_{L^{2}_{t}L^{2}_{x}} \leq \, C\,  T^{\frac{1}{4}} \times T^{\frac{1}{2}-\frac{1}{4}} \left\| \theta \right\|_{L^{4}_{t}L^{2}_{x}}. 
\end{split}
\end{equation*}
\finpv 
\begin{Lemma} We have 
\begin{equation}\label{Bilinear-limit-case}
\left\| \int_{0}^{t} e^{(t-\tau)\Delta}\vu \cdot \vN \theta (\tau,\cdot)d \tau \right\|_{L^{4}_{t}L^{2}_{x}} \leq C \, \| \vu \|_{L^{4}_{t}\dot{H}^{1}_{x}} \| \theta \|_{L^{4}_{t}L^{2}_{x}}, 
\end{equation} with a constant $C>0$ independent of $T$.
\end{Lemma}
\pv  We write
\begin{equation*}
  \begin{split}
 \left\| \int_{0}^{t} e^{(t-\tau)\Delta}\vu \cdot \vN \theta (\tau,\cdot)d \tau \right\|_{L^{4}_{t}L^{2}_{x}} =&\, \left\| \int_{0}^{t} e^{(t-\tau)\Delta} div(\theta \vu) (\tau,\cdot)d \tau \right\|_{L^{4}_{t}L^{2}_{x}}\\
  \leq & \, C\, \left\| \int_{0}^{t} e^{(t-\tau)\Delta} (\theta \vu) (\tau,\cdot)d \tau \right\|_{L^{4}_{t}\dot{H}^{1}_{x}}. 
  \end{split}  
\end{equation*}
Using the third point in Lemma \ref{heat-estimates-2} with $s_1=-1/2$, $s_2=3/2$ and $p=4$, we get 
\begin{equation*}
C\, \left\| \int_{0}^{t} e^{(t-\tau)\Delta} (\theta \vu) (\tau,\cdot)d \tau \right\|_{L^{4}_{t}\dot{H}^{1}_{x}} \leq  C\, \left( \int_{0}^{T} \| \theta \vu \|^{2}_{\dot{H}^{-1/2}} dt \right)^{1/2}.
\end{equation*}
Applying the Hardy-Littlewood-Sobolev inequalities and H\"older inequalities, we write
\begin{equation*}
\begin{split}
 C\, \left( \int_{0}^{T} \| \theta \vu \|^{2}_{\dot{H}^{-1/2}} dt \right)^{1/2}\leq &\, C \, \left( \int_{0}^{T} \| \theta \vu (t,\cdot)\|^{2}_{L^{3/2}} \,  dt \right)^{1/2} \\
\leq & \, C \, \left( \int_{0}^{T} \| \theta(t,\cdot)\|^{2}_{L^2} \, \| \vu (t,\cdot)\|^{2}_{L^{6}} \,  dt \right)^{1/2} \\
\leq & \, C \, \left( \int_{0}^{T} \| \theta(t,\cdot)\|^{2}_{L^2} \, \| \vu (t,\cdot)\|^{2}_{\dot{H}^1} \,  dt \right)^{1/2}\\
\leq & \, C  \, \| \vu \|_{L^{4}_{t}\dot{H}^{1}_{x}} \| \theta \|_{L^{4}_{t}L^{2}_{x}}.
\end{split}
\end{equation*}
\finpv

\noindent
Now, we are able to set a time $T_0$ small and to apply Lemma \ref{Picard} as follows:  first, by estimate (\ref{Linear1-limit-case}) we set $T_0$ small enough to satisfy the first condition in (\ref{Conditions-Picard}). On the other hand, recall that by estimates (\ref{Bilinear-N-S2}) and (\ref{Bilinear-limit-case}) the constant $C_B$ does not depend on $T_0$. Thus, we define $\delta =  \left\| e^{t\Delta} \vu_0 \right\|_{L^{4}_{t}}+\left\|e^{t\Delta} \theta_0 \right\|_{L^{4}_{t}L^{2}_{x}}$, and in estimates (\ref{Control-u0})-(\ref{Control-theta0}) we impose an additional   smallness condition on $T_0$ to get $\delta <\varepsilon$, with $\varepsilon$ small. We thus satisfy the second and the third condition in (\ref{Conditions-Picard}). This way, we obtain a solution $(\vu, \theta) \in L^4([0,T_0], \dot{H}^1(\Rt)) \times L^4([0,T_0], L^2(\Rt) )$ to (\ref{Mild-u})-(\ref{Mild-theta}).  

\medskip

\noindent
This solution also belongs to the space $E_1 \times E_2$. To verify this fact, we shall use the following  result. 
\begin{Lemma}[Lemma $5.10$ of  \cite{BaCheDan}]\label{Lemma-Functional-Space} Let $v \in \mathcal{C}([0,T],\mathcal{S}'(\Rt))$ be a solution of the heat equation
\begin{equation*}
    \partial_t v - \Delta v = g, \quad u(0,\cdot)=u_0.
\end{equation*}
Let $\sigma \in \R$. If $v_0 \in \dot{H}^{\sigma}(\Rt)$ and $g \in L^2([0,T],\dot{H}^{\sigma-1}(\Rt))$ then we have 
\[ v \in L^{\infty}([0,T],\dot{H}^{\sigma}(\Rt))\cap L^2([0,T],\dot{H}^{\sigma+1}(\Rt)). \]
\end{Lemma}

\medskip

\noindent
We start by proving that $\theta \in E_2$. In this lemma we set $\sigma=-\frac{1}{2}$, $v_0=\theta_0\in \dot{H}^{-1/2}(\Rt)$ and $g = \vu \cdot \vN \theta = div(\theta \vu)$, where   we must verify that 
\[ g \in  L^{2}([0,T_0], \dot{H}^{-3/2} (\Rt)). \]
Indeed, since $\vu \in L^{4}_{t}\dot{H}^{1}_x \subset L^{4}_{t}L^{6}_x$ and $\theta \in L^{4}_{t}L^{2}_{x}$ by H\"older inequalities and Hardy-Littlewood-Sobolev inequalities we have $\theta \vu \in L^{2}_{t}L^{3/2}_{x}\subset L^{2}_{t} \dot{H}^{-1/2}_{x}$, hence we get $g=div(\theta \vu)\in L^{2}_{t}\dot{H}^{-3/2}_{x}$. We thus have  $\theta \in E_2$. 

\medskip

\noindent
Now we prove that $\vu \in E_1$. Note that this lemma also holds for vector fields, and by a slight abuse of notation  we set $\sigma=\frac{1}{2}$, $v_0=\vu_0$ and $g= -\P( (\vu \cdot \vN)\vu) + \P(\theta \ve_3)$, where we will verify that 
\[ g\in L^2([0,T_0],\dot{H}^{-1/2}(\Rt)). \]
It is well-known that $\vu \in L^{4}_{t}\dot{H}^{1}_x$ yields $\P( (\vu \cdot \vN)\vu)\in L^{2}_{t}\dot{H}^{-1/2}_x$ (see \cite[Theorem $5.6$]{BaCheDan}). Moreover,
since $\theta \in E_2$ we have $\P(\theta \ve_2) \in L^{\infty}_{t}\dot{H}^{-1/2}_{x} \subset L^{2}_{t}\dot{H}^{-1/2}_{x}$. We thus get $\vu \in L^{\infty}_t \dot{H}^{1/2}_x \cap L^2_t \dot{H}^{3/2}_x$ by Lemma \ref{Lemma-Functional-Space}, and it remains to prove that $\vu \in L^\infty_t L^2_x$. 
Since $(\vu \cdot \vN) \vu=\text{div}(\vu \otimes \vu)\in L^2_t \dot{H}^{-1/2}_x$ we can perform the classical energy estimate in the first equation in the system (\ref{Boussinesq}) to get $\vu \in L^\infty_t L^2_x$. This way, we have $\vu \in E_1$.

\subsubsection{When $s=1/2$ and $1/2<r\leq 1$.}
Here, we consider the spaces
\[ E_1= L^{\infty}([0,T], H^{r}(\Rt))\cap L^{2}([0,T],\dot{H}^{ r+1}(\Rt)) \]
and 
\[ E_2= L^{\infty}([0,T], \dot{H}^{ -\frac{1}{2} }(\Rt))\cap L^{2}([0,T],\dot{H}^{ \frac{1}{2} }(\Rt)),\]
with their usual norms. In order  to construct a solution $(\vu,\theta)\in E_1 \times E_2$ of the system (\ref{Mild-u})-(\ref{Mild-theta}),  we shall follow   some of the  ideas  of the previous case when $s=r=1/2$. First, remark that  we have  the following embeddings:
\[ E_1 \subset  L^{4}([0,T],\dot{H}^{r+1/2}(\Rt)) \cap L^4([0,T], \dot{H}^1(\Rt))=F_1,\]
and 
\[E_2 \subset L^4([0,T],L^2(\Rt))\cap L^{4/(2r-1)}([0,T],\dot{H}^{r-1}(\Rt))=F_2.\]
Indeed, by interpolation inequalities we have $E_1 \subset L^{4}_{t}\dot{H}^{r+1/2}_{x}$. Moreover, we  have $E_1 \subset L^{\infty}_{t}H^{1/2}_{x}\cap L^{2}_{t}\dot{H}^{3/2}_{x}\subset L^{4}_{t}\dot{H}^{1}_{x}$, which yields the first embedding $E_1 \subset F_1$. The second embedding $E_2 \subset F_2$ also follows  from interpolation inequalities: on the one hand, we have $E_2 \subset L^{4}_{t}L^{2}_{x}$, and on the other hand,   since $1/2<r\leq 1$ then $-1/2<r-1\leq 0$ and we have $E_2 \subset L^{4/(2r-1)}_{t}\dot{H}^{r-1}_{x}$.  

\medskip

\noindent
Spaces $F_1$ and $F_2$ are equipped with their standard norms
\[ \| \cdot \|_{F_1}=\| \cdot \|_{L^{4}_{t}\dot{H}^{1}_{x}}+\| \cdot \|_{L^{4}_{t}\dot{H}^{r+1/2}_{x}}, \]
\[ \| \cdot \|_{F_2}=\| \cdot \|_{L^{4}_{t}L^{2}_{x}}+\| \cdot \|_{L^{4/(2r+1)}_{t}\dot{H}^{r-1}_{x}},\]
respectively. Firstly,  these norms  will allow us to control all the terms in equations (\ref{Mild-u})-(\ref{Mild-theta}) and to construct a solution  $(\vu,\theta) \in F_1 \times F_2$. Then, we will verify that this solution belongs to the space $E_1 \times E_2$, and we shall finish this part of the proof with the uniqueness of this solution.

\medskip

\noindent
For the data term  in equation (\ref{Mild-u}), in Lemma \ref{heat-estimates-3}
we set the parameters $(s_1,s_2)=(r,r+1/2)$ and $(s_1,s_2)=(1/2,1)$, hence in both cases we get $p=4$, and we have the control:
\begin{equation}\label{Control-u0-2}
\left\| e^{t\Delta}\vu_0 \right\|_{F_1} \leq \frac{\varepsilon}{2}+(R^2_\varepsilon T)^{1/4}\| \vu_0 \|_{H^r}.  
\end{equation}
Similarly, for the data term in equation (\ref{Mild-theta}), in Lemma \ref{heat-estimates-3} we set now the parameters $(s_1,s_2=-1/2,0)$, hence we get $p=4$, and $(s_1,s_2)=(-1/2,r-1)$, hence we obtain $p=\frac{4}{2r-1}$. We thus have the control:
\begin{equation}\label{Control-theta0-2}
\left\| e^{t\Delta}\theta_0 \right\|_{F_2} \leq \frac{\varepsilon}{2}+\Big((R^2_\varepsilon T)^{1/4}+(R^2_\varepsilon T)^{(2r-1)/4}\Big)\| \theta_0 \|_{\dot{H}^{-1/2}}.  
\end{equation}

\noindent
Then, the rest of the terms in equations (\ref{Mild-u})-(\ref{Mild-theta}) can be estimated as follows. 
\begin{Lemma} We have
\begin{equation}\label{Bilinear-N-S3}
\left\| \int_{0}^{t} e^{(t-\tau)\Delta} \P \left( (\vu \cdot \vN)\vv\right)(\tau,\cdot)d \tau \right\|_{F_1} \leq C \| \vu \|_{F_1} \| \vv \|_{F_1},
\end{equation} 
with a constant $C>0$ independent of $T$.
\end{Lemma}
\pv The term involving the norm $\| \cdot \|_{L^{4}_{t}\dot{H}^{1}_{x}}$ was estimated in (\ref{Bilinear-N-S2}). Then, for the  norm $\| \cdot \|_{L^{4}_{t}\dot{H}^{r+1/2}}$, by the third point of Lemma \ref{heat-estimates-2}  we get 
\[\left\| \int_{0}^{t} e^{(t-\tau)\Delta} \P \left( (\vu \cdot \vN)\vv\right)(\tau,\cdot)d \tau \right\|_{L^{4}_{t}\dot{H}^{r+1/2}} \leq  C\, \|  (\vu \cdot \vN)\vv \|_{L^{2}_{t}\dot{H}^{r-1}_{x}},\]
and we write
\begin{equation}\label{Estim-tech-01}
\begin{split}
C\, \|  (\vu \cdot \vN)\vv \|_{L^{2}_{t}\dot{H}^{r-1}_{x}} \leq  &\, \, C\, \| div(\vv \otimes \vu) \|_{L^{2}_{t}\dot{H}^{r-1}_{x}}\\
\leq &\, C\, \| \vv \otimes \vu \|_{L^{2}_{t}\dot{H}^{r}_{x}}.
\end{split}
\end{equation}
By product laws in homogeneous Sobolev spaces (with the relationship $r=(r+1/2)+1 - 3/2$), and by H\"older inequalities in the time variable (with $1/2=1/4+1/4$) we obtain
\begin{equation}\label{Estim-tech-02}
\begin{split}
C\, \| \vv \otimes \vu \|_{L^{2}_{t}\dot{H}^{r}_{x}} \leq &\,  C \| \vv \|_{L^{4}_{t}\dot{H^{r+1/2}_{x}}}\| \vu \|_{L^{4}_{t}\dot{H}^{1}_{x}}+C \| \vu \|_{L^{4}_{t}\dot{H^{r+1/2}_{x}}}\| \vv \|_{L^{4}_{t}\dot{H}^{1}_{x}}\\
\leq &\, C\, \| \vu \|_{F_1}\, \| \vv \|_{F_1}, 
\end{split}
\end{equation}
hence estimate (\ref{Bilinear-N-S3}) follows.  \finpv 

\begin{Lemma} It holds
\begin{equation}\label{Linear-2}
\left\| \int_{0}^{t}e^{(t-\tau)\Delta} \P \Big( \theta \ve_3 \Big)(\tau,\cdot)\right\|_{F_1} \leq C\max\left( T^{1/2},T^{(3-2r)/4}\right)  \| \theta \|_{F_2},
\end{equation}
where, since $1/2<r\leq 1$ we have $(3-2r)/4>0$.
\end{Lemma}
\pv  The term concerning the norm $\| \cdot \|_{L^{4}_{t}\dot{H}^{1}_{x}}$ was estimated in (\ref{Linear1-limit-case}), so it remains to study the term concerning the norm $\| \cdot \|_{L^{4}_{t}\dot{H}^{r+1/2}_{x}}$. By the third point of Lemma \ref{heat-estimates-2}, and by the fact that $4/(2r-1)>2$ (recall that  $1/2<r\leq 1$), we write
\begin{equation*}
\begin{split}
\left\| \int_{0}^{t}e^{(t-\tau)\Delta} \P \Big( \theta \ve_3 \Big)(\tau,\cdot)\right\|_{L^{4}_{t}\dot{H}^{r+1/2}_{x}} \leq &\,  C \, \| \theta \|_{L^{2}_{t}\dot{H}^{r-1}_{x}} \\
\leq &\,  C\,T^{1/2 - (2r-1)/4}\| \theta \|_{L^{4/(2r-1)}_{t}\dot{H}^{r-1}_{x}},
\end{split}
\end{equation*}
where the expression $1/2 - (2r-1)/4$ computes down as $(3-2r)/4$. We thus obtain the wished estimate (\ref{Linear-2}).  \finpv 

\begin{Lemma} We have
\begin{equation}\label{Bilinear-2}
\left\| \int_{0}^{t} e^{(t-\tau)\Delta}\vu \cdot \vN \theta (\tau,\cdot)d \tau \right\|_{F_2} \leq C(1+T^{(2r-1)/4}) \, \| \vu \|_{F_1} \| \theta \|_{F_2}. 
\end{equation}
\end{Lemma}
\pv The term involving the norm $\| \cdot \|_{L^{4}_{t}L^{2}_{x}}$ was already treated in estimate (\ref{Bilinear-limit-case}), and must  study   the term involving the norm $\| \cdot \|_{L^{4/(2r-1)}_{t}\dot{H}^{r-1}_{x}}$. For this, remark first that in the case $r=1$ we have $L^{4/(2r-1)}_{t}\dot{H}^{r-1}_{x}=L^{4}_{t}L^{2}_{x}$, so this case is done. We thus focus on the case $1/2<r<1$. We start by writing
\begin{equation*}
\begin{split}
&\, \left\| \int_{0}^{t} e^{(t-\tau)\Delta}\vu \cdot \vN \theta (\tau,\cdot)d \tau \right\|_{L^{4/(2r-1)}_{t}\dot{H}^{r-1}_{x}} \\
\leq &\, C\, T^{(2r-1)/4} \left\| \int_{0}^{t} e^{(t-\tau)\Delta}\vu \cdot \vN \theta (\tau,\cdot)d \tau \right\|_{L^{\infty}_{t}\dot{H}^{r-1}_{x}} \\
\leq &\, C\, T^{(2r-1)/4} \left\| \int_{0}^{t} e^{(t-\tau)\Delta} div \left( (-\Delta)^{(r-1)/2}(\theta \vu)\right) (\tau,\cdot)d \tau \right\|_{L^{\infty}_{t}L^{2}_{x}}. 
\end{split}
\end{equation*}
We apply the first point of Lemma \ref{heat-estimates-2} to obtain
\begin{equation*}
\begin{split}
&\, C\, T^{(2r-1)/4} \left\| \int_{0}^{t} e^{(t-\tau)\Delta} div \left( (-\Delta)^{(r-1)/2}(\theta \vu)\right) (\tau,\cdot)d \tau \right\|_{L^{\infty}_{t}L^{2}_{x}}\\
\leq &\, C\, T^{(2r-1)/4} \| (-\Delta)^{(r-1)/2}(\theta \vu) \|_{L^{2}_{t}L^{2}_{x}}= C\, T^{(2r-1)/4} \| \theta \vu \|_{L^{2}_{t}\dot{H}^{r-1}_{x}}.
\end{split}
\end{equation*}
Observe that we $-1/2<r-1<0$ (since $1/2<r<1$). We thus apply Hardy-Littlewood-Sobolev inequalities and for $p=6/(5-2r)$ (which verifies $3/2<p<2$) we write
\begin{equation*}
C\, T^{(2r-1)/4} \| \theta \vu \|_{L^{2}_{t}H^{r-1}_{x}} \leq C\, T^{(2r-1)/4} \| \theta \vu \|_{L^{2}_{t}L^{p}_{x}}.  
\end{equation*} 
Then, we use H\"older inequalities (with $1/p=1/2+1/q$ and $q=3/(1-r)$) to obtain
\begin{equation*}
C\, T^{(2r-1)/4} \| \theta \vu \|_{L^{2}_{t}L^{p}_{x}} \leq C\, T^{(2r-1)/4}  \left( \int_{0}^{T} \| \theta (\tau,\cdot)\|^{2}_{L^2}\| \vu(\tau,\cdot)\|^{2}_{L^q}d \tau \right)^{1/2}.
\end{equation*}
In the last term, remark that by Hardy-Littlewood-Sobolev inequalities we have the continuous embedding $\dot{H}^{r+1/2}(\Rt) \subset L^q(\Rt)$. Moreover, by H\"older inequalities in the time variable (with $1/2=1/4+1/4$), we have  
\begin{equation*}
\begin{split}
&\, C\, T^{(2r-1)/4}  \left( \int_{0}^{T} \| \theta (\tau,\cdot)\|^{2}_{L^2}\| \vu(\tau,\cdot)\|^{2}_{L^q}d \tau \right)^{1/2} \\
\leq &\,  C\, T^{(2r-1)/4}  \left( \int_{0}^{T} \| \theta (\tau,\cdot)\|^{2}_{L^2}\| \vu(\tau,\cdot)\|^{2}_{\dot{H}^{r+1/2}}d \tau \right)^{1/2}\\
\leq &\, C\, T^{(2r-1)/4}\, \| \theta \|_{L^{4}_{t}L^{2}_{x}}\| \| \vu \|_{L^{4}_{t}\dot{H}^{r+1/2}_{x}}\\
\leq &\, C\, T^{(2r-1)/4}\, \| \theta \|_{F_2}\| \| \vu \|_{F_1}.
\end{split}
\end{equation*}
We thus obtain the desired estimate (\ref{Bilinear-2}). \finpv 

\noindent
With estimates (\ref{Control-u0-2}), (\ref{Control-theta0-2}), (\ref{Bilinear-N-S3}), (\ref{Linear-2}) and (\ref{Bilinear-2}) at hand, we proceed as in the previous case (when $r=s=1/2$) to obtain a solution $(\vu, \theta)\in F_1 \times F_2$ of the system (\ref{Mild-u})-(\ref{Mild-theta}) for a time $T_0>0$ small. 

\medskip

\noindent
Now, by Lemma \ref{Lemma-Functional-Space} and the information $(\vu, \theta)\in F_1 \times F_2$ we have $(\vu, \theta)\in E_1 \times E_2$. Indeed, the fact that $\theta \in E_2$ was already proven in Section \ref{Sec:s-r-1/2} (below Lemma \ref{Lemma-Functional-Space}), so it remains to prove that $\vu \in E_1$.  For this, we must  prove that $(\vu \cdot \vN) \vu + \theta \ve_3 \in L^{2}_{t}\dot{H}^{r-1}_{x}$. For the first term, by estimates (\ref{Estim-tech-01}) and (\ref{Estim-tech-02}) we directly obtain $(\vu \cdot \vN) \vu \in L^{2}_{t}\dot{H}^{r-1}_{x}$. For the second term, since $\theta \in F_2$ we have $\theta \in L^{4/(2r-1)}_{t}\dot{H}^{r-1}_{x}$. Moreover, since $r\leq 1$ in particular we have $r<3/2$ hence $4/(2r-1)>2$. We thus get $\theta \in L^{2}_{t}\dot{H}^{r-1}_{x}$, and by Lemma \ref{Lemma-Functional-Space} we have $\vu \in L^{\infty}_t \dot{H}^r_x \cap L^2_t \dot{H}^{r+1}_x$. Moreover, following the same ideas  at the end of Section \ref{Sec:s-r-1/2}, we also have  $\vu \in L^\infty_t L^2_x$ and we get $\vu\in E_1$. 

\medskip

\noindent
Theorem \ref{Th1} is now proven. \finpv

\section{Proof of Proposition \ref{Prop-Uniqueness}}\label{Sec:Proof-Prop-Uniq}
As in Section \ref{Sec:s-r-1/2}, we denote 
\[ E_1= L^{\infty}([0,T], H^{ \frac{1}{2} }(\Rt))\cap L^{2}([0,T],\dot{H}^{ \frac{3}{2} }(\Rt)), \]
and 
\[ E_2= L^{\infty}([0,T], \dot{H}^{ -\frac{1}{2} }(\Rt))\cap L^{2}([0,T],\dot{H}^{ \frac{1}{2}}\cap \dot{W}^{1,3}(\Rt)).\]
Moreover,  we define $\vv=\vu_1 - \vu_2$ and $\eta=\theta_1 - \theta_2$. Then,  the couple $(\vv,\eta)$ solves the following system:
\begin{equation*}
\begin{cases}\vspace{2mm}
    \partial_t \vv -\Delta \vv + \P\big( (\vu_1 \cdot \vN) \vv + (\vv \cdot \vN)\vu_2 \big)-\P(\eta \ve_3)=0, \ \ div(\vv)=0, \\ \vspace{2mm}
    \partial_t \eta -\Delta \eta + \vu_1 \cdot \vN \eta + \vv \cdot \vN \theta_2 =0,\\    
    \vv(0,\cdot)=0, \quad \eta(0,\cdot)=0.
 \end{cases}   
\end{equation*}

\medskip

\noindent
Since $\vv \in E_1$ and $\eta \in E_2$, in the first equation above we can perform an energy estimate in the $\dot{H}^{1/2}-$inner product  $\langle f,g\rangle_{\dot{H}^{1/2}}=\int_{\Rt} |\xi| \widehat{f}(\xi) \Bar{\widehat{g}}(\xi) d \xi$, while in the second equation above we perform an energy estimate in the $\dot{H}^{-1/2}-$inner product  $\langle f,g\rangle_{\dot{H}^{-1/2}}=\int_{\Rt} |\xi|^{-1} \widehat{f}(\xi) \Bar{\widehat{g}}(\xi) d \xi$. Then, for $0<t\leq T_0$ we obtain 
\begin{equation*}
\begin{split}
&\, \frac{d}{dt} \|\vv(t,\cdot)\|^{2}_{\dot{H}^{1/2}}+2 \| \vN \otimes \vv(t,\cdot)\|^{2}_{\dot{H}^{1/2}}  \\
= &\, -2 \Big\langle(\vu_1 \cdot \vN) \vv + (\vv \cdot \vN)\vu_2, \vv\Big\rangle_{\dot{H} ^{1/2}}  +2   \Big\langle \eta \ve_3, \vv\Big\rangle_{\dot{H} ^{1/2}},
\end{split}
\end{equation*}
\begin{equation*}
\begin{split}
&\, \frac{d}{dt}\| \eta(t,\cdot)\|^{2}_{\dot{H}^{-1/2}}+2\| \vN \eta(t,\cdot)\|^2_{\dot{H}^{-1/2}}\\
=&\, -2\Big\langle  \vu_1 \cdot \vN \eta + \vv \cdot \vN \theta_2, \eta\Big\rangle_{\dot{H} ^{-1/2}}. 
\end{split}
\end{equation*}
To simplify our writing, we shall denote 
\begin{equation*}
 \mathcal{E}_1(t)=\|\vv(t,\cdot)\|^{2}_{\dot{H}^{1/2}} + \|\eta(t,\cdot)\|^{2}_{\dot{H}^{-1/2}},    
\end{equation*}
and
\begin{equation*}
 \mathcal{E}_2(t)= \| \vN \otimes \vv(t,\cdot)\|^{2}_{\dot{H}^{1/2}}+\| \vN \eta(t,\cdot)\|^2_{\dot{H}^{-1/2}}.   
\end{equation*}
We thus have
\begin{equation}\label{Iden1}
  \begin{split}
&\frac{d}{dt} \mathcal{E}_1(t)+2\mathcal{E}_2(t)\\
=&\, -2 \Big\langle(\vu_1 \cdot \vN) \vv + (\vv \cdot \vN)\vu_2, \vv\Big\rangle_{\dot{H} ^{1/2}}  +2   \Big\langle \eta \ve_3, \vv\Big\rangle_{\dot{H} ^{1/2}}\\
&\,-2\Big\langle  \vu_1 \cdot \vN \eta + \vv \cdot \vN \theta_2, \eta\Big\rangle_{\dot{H} ^{-1/2}},
  \end{split}  
\end{equation}
where we must estimate each term on the right-hand side. 
\begin{Lemma} We have 
\begin{equation}\label{Term1}
\begin{split}
&\, 2\left| \Big\langle (\vu_1 \cdot \vN) \vv + (\vv \cdot \vN)\vu_2, \vv\Big\rangle_{\dot{H} ^{1/2}} \right|\\
\leq &\, C\, \mathcal{E}_1(t)  \Big( \| \vu_1(t,\cdot)\|^{4}_{\dot{H}^{1}} + \|\vu_2(t,\cdot)\|^{4}_{\dot{H}^{1}}  \Big) + \frac{1}{3} \mathcal{E}_2(t). 
\end{split}
\end{equation}
\end{Lemma}
\pv  By \cite[Lemma $5.12$]{BaCheDan}, the interpolation inequalities and the discrete Young inequalities (with $1=1/4+3/4$),  we have
\begin{equation*}
\begin{split}
&\, 2\left| \Big\langle (\vu_1 \cdot \vN)\vv, \vv\Big\rangle_{\dot{H} ^{1/2}} \right| \\
\leq &\,  C\,  \| \vu_1(t,\cdot)\|_{\dot{H}^1} \| \vv(t,\cdot)\|_{\dot{H}^{1}}\| \vN \otimes \vv (t,\cdot)\|_{\dot{H}^{1/2}} \\
\leq &\, C\, \| \vu_1(t,\cdot)\|_{\dot{H}^{1}}\| \vv(t,\cdot)\|^{1/2}_{\dot{H}^{1/2}} \| \vv(t,\cdot)\|^{1/2}_{\dot{H}^{3/2}}\, \| \vN \otimes \vv (t,\cdot)\|_{\dot{H}^{1/2}}\\
\leq &\, C\, \| \vu_1(t,\cdot)\|_{\dot{H}^{1}} \| \vv(t,\cdot)\|^{1/2}_{\dot{H}^{1}}\| \vN \otimes \vv (t,\cdot)\|^{3/2}_{\dot{H}^{1/2}}\\
\leq &\, C\, \| \vu_1(t,\cdot)\|^{4}_{\dot{H}^{1}} \| \vv(t,\cdot)\|^{2}_{\dot{H}^{1/2}}+ \frac{1}{6}\| \vN \otimes \vv (t,\cdot)\|^{2}_{\dot{H}^{1/2}}\\
\leq &\, C \mathcal{E}_1(t)\| \vu_1(t,\cdot)\|^{4}_{\dot{H}^{1}}+\frac{1}{6}\mathcal{E}_2(t).
\end{split}
\end{equation*}
By the same arguments, we can write
\begin{equation*}
    \begin{split}
&\, 2\left| \Big\langle (\vv \cdot \vN)\vu_2, \vv\Big\rangle_{\dot{H} ^{1/2}} \right| \\
\leq &\,  C\,  \| \vv(t,\cdot)\|_{\dot{H}^1} \| \vu_2(t,\cdot)\|_{\dot{H}^{1}}\| \vN \otimes \vv (t,\cdot)\|_{\dot{H}^{1/2}} \\
\leq &\, C\, \| \vv(t,\cdot)\|^{1/2}_{\dot{H}^{1/2}}\| \vv(t,\cdot)\|^{1/2}_{\dot{H}^{3/2}}\, \| \vu_2(t,\cdot)\|_{\dot{H}^{1}}\| \vN \otimes \vv (t,\cdot)\|_{\dot{H}^{1/2}}\\
\leq &\, C\, \| \vv(t,\cdot)\|^{1/2}_{\dot{H}^{1/2}} \| \vu_2(t,\cdot)\|_{\dot{H}^{1}}\| \vN \otimes \vv (t,\cdot)\|^{3/2}_{\dot{H}^{1/2}}\\
\leq &\, C\, \| \vv(t,\cdot)\|^{2}_{\dot{H}^{1/2}} \| \vu_2(t,\cdot)\|^{4}_{\dot{H}^{1}} + \frac{1}{6} \| \vN \otimes \vv (t,\cdot)\|^{2}_{\dot{H}^{1/2}}\\
\leq &\, C \mathcal{E}_1(t)\| \vu_2(t,\cdot)\|^{4}_{\dot{H}^{1}}+\frac{1}{6}\mathcal{E}_2(t).
    \end{split}
\end{equation*}
Gathering these estimates, we obtain (\ref{Term1}).  \finpv 
\begin{Lemma} It holds
\begin{equation}\label{Term2}
\left| \Big\langle \eta \ve_3, \vv\Big\rangle_{\dot{H} ^{1/2}} \right| \leq C \mathcal{E}_1(t)+\frac{1}{3}\mathcal{E}_2(t).
\end{equation}    
\end{Lemma}
\pv  We write
\begin{equation*}
    \begin{split}
    2 \left| \Big\langle \eta \ve_3, \vv\Big\rangle_{\dot{H} ^{1/2}} \right| \leq &\, 2 \int_{\Rt} |\widehat{\eta}(t,\xi)|  \, |\xi| |\widehat{\vv}(t,\xi)| d \xi \\
    \leq &\, 2 \| \eta(t,\cdot)\|_{L^2}\, \| \vv(t,\cdot)\|_{\dot{H}^{1}}\\
    \leq &\, C\| \eta(t,\cdot)\|_{L^2} \| \vv(t,\cdot)\|^{1/2}_{\dot{H}^{1/2}}\, \| \vN \otimes \vv(t,\cdot)\|^{1/2}_{\dot{H}^{1/2}}\\
    =&\, C \| \vv(t,\cdot)\|^{1/2}_{\dot{H}^{1/2}}\, \| \eta(t,\cdot)\|_{L^2}\| \vN \otimes \vv(t,\cdot)\|^{1/2}_{\dot{H}^{1/2}}.
    \end{split}
\end{equation*}
Then, we apply the discrete Young inequalities (first with $1=1/4+3/4$ and thereafter with $1=1/3+2/3$) to get
\begin{equation*}
    \begin{split}
 &\,  C \| \vv(t,\cdot)\|^{1/2}_{\dot{H}^{1/2}}\, \| \eta(t,\cdot)\|_{L^2}\| \vN \otimes \vv(t,\cdot)\|^{1/2}_{\dot{H}^{1/2}} \\
   \leq &\,C\,  \| \vv(t,\cdot)\|^2_{\dot{H}^{1/2}}+ \|\eta(t,\cdot)\|^{4/3}_{L^2}\, \| \vN\otimes \vv(t,\cdot)\|^{2/3}_{\dot{H}^{1/2}}\\
   \leq & \, C\,  \| \vv(t,\cdot)\|^2_{\dot{H}^{1/2}} + C\, \|\eta(t,\cdot)\|^2_{L^2}+\frac{1}{6}\| \vN\otimes \vv(t,\cdot)\|^{2}_{\dot{H}^{1/2}}\\
   \leq &\, C\, \mathcal{E}_1(t)+C\, \|\eta(t,\cdot)\|^2_{L^2}+\frac{1}{6}\mathcal{E}_2(t).
    \end{split}
\end{equation*}
Moreover, by interpolation inequalities and using again the discrete Young inequalities (with $1=1/2+1/2$), the term $\|\eta(t,\cdot)\|^2_{L^2}$ is estimated as:
\begin{equation*}
\begin{split}
C\, \|\eta(t,\cdot)\|^2_{L^2} \leq &\,  C\, \| \eta(t,\cdot)\|_{\dot{H}^{-1/2}}\| \eta(t,\cdot)\|_{\dot{H}^{1/2}} \\
\leq &\, C \, \| \eta(t,\cdot)\|_{\dot{H}^{-1/2}}\| \vN \eta(t,\cdot)\|_{\dot{H}^{-1/2}}\\
\leq &\, C \, \| \eta(t,\cdot)\|_{\dot{H}^{-1/2}}^2 + \frac{1}{6}\| \vN \eta(t,\cdot)\|^2_{\dot{H}^{-1/2}}\\
\leq &\, C\, \mathcal{E}_1(t)+\frac{1}{6}\mathcal{E}_2(t).
\end{split}
\end{equation*}
Gathering these estimates, we obtain (\ref{Term2}). \finpv 

\begin{Lemma}  We have
 \begin{equation}\label{Term3}
\begin{split}
 &\, 2\left|\Big\langle  \vu_1 \cdot \vN \eta + \vv \cdot \vN \theta_2, \eta\Big\rangle_{\dot{H} ^{-1/2}}\right| \\
\leq &\, C\, \mathcal{E}_1(t) \Big( \| \vu_1(t,\cdot)\|^4_{\dot{H}^1}+\|\theta_2(t,\cdot)\|^2_{\dot{W}^{1,3}}+1\Big)+ \frac{1}{3} \mathcal{E}_2(t).
\end{split}
\end{equation}   
\end{Lemma}
\pv To study the first term involving the expression $\vu_1 \cdot \vN \eta$, by the identity $\vu_1 \cdot \vN \eta= div(\eta \vu_1)$ we write 
\begin{equation*}
 \begin{split}
 2 \left| \Big\langle  div(\eta \vu_1), \eta \Big\rangle_{\dot{H}^{-1/2}} \right| \leq &\, C\int_{\Rt}|\xi|^{-1}|\xi| | \widehat{\eta \vu_1} |\, |\widehat{\eta}| d \xi \\
 \leq &\, C \,  \int_{\Rt}| \widehat{\eta \vu_1} |\, |\widehat{\eta}| d \xi\\
 \leq & C\,\| \eta \vu_1(t,\cdot)\|_{L^2}\, \| \eta(t,\cdot)\|_{L^2}.
 \end{split}
\end{equation*}
Then, applying H\"older inequalities and  Hardy-Littlewood-Sobolev inequalities, we get
\begin{equation*}
   \begin{split}
   C\,\| \eta \vu_1(t,\cdot)\|_{L^2}\, \| \eta(t,\cdot)\|_{L^2} \leq &\,C\, \|\eta(t,\cdot)\|_{L^3}\| \vu_1(t,\cdot)\|_{L^6}\|\eta(t,\cdot)\|_{L^2}\\
 \leq &\,  C\, \|\eta(t,\cdot)\|_{\dot{H}^{1/2}}\| \vu_1(t,\cdot)\|_{\dot{H}^1}\|\eta(t,\cdot)\|_{L^2}. 
   \end{split} 
\end{equation*}
Thereafter, in the last we apply interpolation inequalities to write
\begin{equation*}
 \begin{split}
   &\,  C\, \|\eta(t,\cdot)\|_{\dot{H}^{1/2}}\| \vu_1(t,\cdot)\|_{\dot{H}^1}\|\eta(t,\cdot)\|_{L^2} \\
   \leq &\, C\, \|\eta(t,\cdot)\|_{\dot{H}^{1/2}}\| \vu_1(t,\cdot)\|_{\dot{H}^1} \|\eta(t,\cdot)\|^{1/2}_{\dot{H}^{-1/2}}\| \eta(t,\cdot)\|^{1/2}_{\dot{H}^{1/2}}\\
   \leq & \, C\,\|\eta(t,\cdot)\|^{1/2}_{\dot{H}^{-1/2}}\|\vu_1(t,\cdot)\|_{\dot{H}^1} \|  \eta(t,\cdot)\|^{3/2}_{\dot{H}^{1/2}}\\
   \leq & \, C\,\|\eta(t,\cdot)\|^{1/2}_{\dot{H}^{-1/2}}\|\vu_1(t,\cdot)\|_{\dot{H}^1} \| \vN \eta(t,\cdot)\|^{3/2}_{\dot{H}^{-1/2}}.
 \end{split}   
\end{equation*}
We use the discrete Young inequalities (with $1=1/4+3/4$) to get
\begin{equation*}
\begin{split}
&\, C\,\|\eta(t,\cdot)\|^{1/2}_{\dot{H}^{-1/2}}\|\vu_1(t,\cdot)\|_{\dot{H}^1} \| \vN \eta(t,\cdot)\|^{3/2}_{\dot{H}^{-1/2}} \\
\leq &\, C\,\|\eta(t,\cdot)\|^{2}_{\dot{H}^{-1/2}}\|\vu_1(t,\cdot)\|^{4}_{\dot{H}^1}+\frac{1}{3} \| \vN \eta(t,\cdot)\|^{2}_{\dot{H}^{-1/2}}\\
\leq &\, C\, \mathcal{E}_1(t)\,\|\vu_1(t,\cdot)\|^{4}_{\dot{H}^1}+\frac{1}{3}\mathcal{E}_2(t).
 \end{split}
\end{equation*}

\medskip

\noindent
Now, we study the second term involving the expression $\vv \cdot \vN \theta_2$. Using Hardy-Littlewood-Sobolev inequalities and H\"older inequalities (with $2/3=1/3+1/3$) we write
\begin{equation*}
    \begin{split}
\, 2\left| \Big\langle \vv \cdot \vN \theta_2, \eta \Big\rangle_{\dot{H}^{-1/2}} \right|
\leq &\, C \| \vv \cdot \vN \theta_2 \|_{\dot{H}^{-1/2}}\, \| \eta \|_{\dot{H}^{-1/2}}\\
\leq &\, C  \| \vv \cdot \vN \theta_2 \|_{L^{3/2}}\, \| \eta \|_{\dot{H}^{-1/2}}\\
\leq &\, C \| \vv \|_{L^3}\| \vN \theta_2 \|_{L^3}\, \| \eta \|_{\dot{H}^{-1/2}}\\
\leq &\, C \| \vv \|_{\dot{H}^{1/2}}\| \theta_2 \|_{\dot{W}^{1,3}}\,\| \eta \|_{\dot{H}^{-1/2}}.
    \end{split}
\end{equation*}
Then, we apply the discrete Young inequalities to obtain 
\begin{equation*}
\begin{split}
 C \| \vv \|_{\dot{H}^{1/2}}\| \theta_2 \|_{\dot{W}^{1,3}}\,\| \eta \|_{\dot{H}^{-1/2}} \leq &\,  C \| \vv \|^2_{\dot{H}^{1/2}}\| \theta_2 \|^2_{\dot{W}^{1,3}} + \| \eta \|^2_{\dot{H}^{-1/2}}\\
 \leq &\, C \, \mathcal{E}_1(t)\big(\| \theta_2 \|^2_{\dot{W}^{1,3}}+1\big).  
 \end{split}
 \end{equation*}
Gathering these estimates we have the wished inequality (\ref{Term3}). \finpv 

\noindent
With estimates (\ref{Term1}), (\ref{Term2}) and (\ref{Term3}) at our disposal, we get back to identity (\ref{Iden1}) to obtain
\begin{equation*}
\frac{d}{dt}\mathcal{E}_1(t)+\mathcal{E}_2(t) \leq  C \Big(\| \vu_1(t,\cdot)\|^{4}_{\dot{H}^1}+\|\vu_2(t,\cdot)\|^{4}_{\dot{H}^1}+\|\theta_2(t,\cdot)\|^2_{\dot{W}^{1,3}} +1\Big)\mathcal{E}_1(t).
\end{equation*}
Moreover, we define the quantity
\[ \mathcal{N}(t)= \mathcal{E}_1(t)+\int_{0}^{t}\mathcal{E}_2(\tau) d\tau, \]
and by the last inequality we write
\begin{equation*}
    \frac{d}{dt} \mathcal{N}(t) \leq C \Big( \| \vu_1(t,\cdot)\|^{4}_{\dot{H}^1}+\|\vu_2(t,\cdot)\|^{4}_{\dot{H}^1}+\|\theta_2(t,\cdot)\|^2_{\dot{W}^{1,3}}+1\Big) \mathcal{N}(t).
\end{equation*}
We apply the Gr\"onwall inequality to obtain
\begin{equation*}
\begin{split}
&\,\mathcal{N}(t)\\
\leq &\, \mathcal{N}(0)\, \text{exp}\left(C\int_{0}^{t} \Big( \| \vu_1(\tau,\cdot)\|^{4}_{\dot{H}^1}+\|\vu_2(\tau,\cdot)\|^{4}_{\dot{H}^1}+\|\theta_2(\tau,\cdot)\|^2_{\dot{W}^{1,3}}+1\Big)d\tau \right)\\
\leq & \,\mathcal{N}(0)\,\text{exp}\left(C\Big( \| \vu_1\|^{4}_{L^{4}_{t}\dot{H}^1}+\|\vu_2\|^{4}_{L^{4}_{t}\dot{H}^1}+\|\theta_2\|^{2}_{L^2_t \dot{W}^{1,3}_x}+T_0 \Big)\right).
\end{split}
\end{equation*}
Here, recall that $\vu_1,\vu_2 \in E_1 \subset  L^{4}_{t}\dot{H}^{1}_{x}$ and $\theta_2 \in L^2_t \dot{W}^{1,3}_x$, therefore the exponential term above  is well-defined. 
Finally, since $\mathcal{N}(0)=0$ we obtain the wished identities $\vv=0$ and $\eta=0$.  Proposition \ref{Prop-Uniqueness} is proven. \finpv



\begin{thebibliography}{40}
\bibitem{BaCheDan} H. Bahouri, J.-Y. Chemin and  R. Danchin,
\textit{Fourier Analysis and Nonlinear Partial Differential Equations},
Springer-Verlag Berlin Heidelberg 2011, https://doi.org/10.1007/978-3-642-16830-7.
    \bibitem{CaDi} J.-R. Cannon and E. Dibenedetto,
    \textit{The initial value problem for the Boussinesq equations with data in $L^p$}, Approximation methods for Navier-Stokes problems (Proc. Sympos., Univ. Paderborn, Paderborn, 1979), Lecture Notes in Math. 771, Springer, Berlin, 1980, pp. 129–144.

    \bibitem{BranHe} L. Brandolese and J. He, 
    \textit{Uniqueness theorems for the Boussinesq system},
Tohoku Math. J. (2) 72 no.2 (2020), 283--297

    \bibitem{ChYa} D. Chamorro and M. Yangari,
    \textit{Some existence and regularity results for a non-local transport-diffusion equation with fractional derivatives in time and space}, Preprint arXiv:2203.13101 (2022)

    \bibitem{DaPa} R. Danchin and M. Paicu,
    \textit{Les théorèmes de Leray et de Fujita-Kato pour le système de Boussinesq partiellement visqueux}, Bulletin de la Société Mathématique de France, Volume 136 (2008) no. 2, pp. 261-309. doi : 10.24033/bsmf.2557. http://www.numdam.org/articles/10.24033/bsmf.2557/

    \bibitem{DeCu} C. Deng and S. Cui,
    \textit{Well-posedness of the viscous boussinesq system in besov spaces of negative regular index $s = -1$}, J. Math. Phys. 53, 073101 (2012) doi : https://doi.org/10.1063/1.4732521
    

    \bibitem{PFe} P.-G. Fernández-Dalgo
    \textit{Micropolar fluids starting from initial angular velocities with negative Sobolev regularity}, Preprint 	arXiv:2401.16554 (2024)
    
    \bibitem{FuKa} H. Fujita and T. Kato,
    \textit{On the Navier-Stokes initial value problem. I}, Archive for Rational Mechanics and Analysis volume 16, pages 269–315 (1964)

    \bibitem{YKa} Y. Kagei, 
    \textit{On weak solutions of nonstationary Boussinesq equations}, Differential and Integral Equations, 6 (1993), 587–611.

    \bibitem{KLN} K. Kang, J. Lee and D. Duong Nguyen, 
    \textit{Global well-posedness and stability of the 2D Boussinesq equations with partial dissipation near a hydrostatic equilibrium}, Preprint 	arXiv:2306.08286 (2023)

    \bibitem{LPZ} M.-J. Lai, R. Pan and K. Zhao, 
    \textit{Initial Boundary Value Problem for Two-Dimensional Viscous Boussinesq Equations}, Arch Rational Mech Anal 199, 739–760 (2011).

    \bibitem{PLe} P.-G. Lemari\'e-Rieusset. \emph{The Navier-Stokes Problem in the 21st Century}, Chapman \& Hall/CRC, (2016).


    \bibitem{Mo} H. Morimoto, 
    \textit{Non-stationary Boussinesq equations}, J. Fac. Sci. Univ. Tokyo Sect. IA Math 39 (1992), 61–75

    \bibitem{Pedlosky} J. Pedlosky,
    \textit{Geophysical fluid dynamics}, Springer, (1987).

    \bibitem{Salmon} R. Salmon,
    \textit{Lectures on geophysical fluid dynamics}, Oxford University Press, (1998).

    \bibitem{HoLi} Y.-Z. Thomas Hou and C.M Li,
    \textit{Global well-posedness of the viscous Boussinesq equations}, Discr. Cont. Dynam. Sys. 12, 1–12, (2005).
\end{thebibliography}
\end{document}